\documentclass[mathematics,article,accept,moreauthors,pdftex]{Definitions/mdpi} 
\usepackage{comment}
\usepackage{caption}
\usepackage{subfig}
\usepackage{graphicx}
\usepackage{enumitem}
\usepackage{array}
\usepackage{smartdiagram}
\usepackage{amssymb}
\usepackage{color}
\colorlet{grey1}{black!10!white!90!}
\colorlet{grey2}{black!20!white!80!}
\colorlet{grey3}{black!30!white!70!}
\colorlet{grey4}{black!40!white!60!}
\colorlet{grey5}{black!50!white!50!}
\colorlet{grey6}{black!60!white!40!}
\colorlet{grey7}{black!70!white!30!}
\definecolor{grey5}{rgb}{0.52, 0.52, 0.51}
\definecolor{grey8}{rgb}{0.21, 0.27, 0.31}
\colorlet{grey9}{black!90!white!10!}
\usepackage{pict2e,epic,eepic}
\usepackage{pstricks}
\usepackage{babel} 
\usepackage{comment} 
\usetikzlibrary{babel}
\usetikzlibrary{shapes.geometric}
\usetikzlibrary{shapes.arrows}
\usetikzlibrary{arrows,calc,positioning,fit}
\usepackage{pgf,pgfplots}
\pgfplotsset{compat=1.14}
\usepackage{tikz,pgfkeys}
\usepackage{algorithm}
\usepackage{algorithmic}

\firstpage{1} 
\pubvolume{1}
\issuenum{1}
\articlenumber{0}
\pubyear{2022}
\copyrightyear{2022}
\externaleditor{Academic Editor: Sławomir T. Wierzchoń}
\datereceived{29 September 2022} 
\dateaccepted{11 November 2022} 
\datepublished{} 
\hreflink{https://doi.org/} 

\Title{From Fuzzy Information to Community Detection: An Approach to Social Networks Analysis with Soft Information}

\TitleCitation{From Fuzzy Information to Community Detection: An Approach to Social Networks Analysis with Soft Information}

\Author{Inmaculada 
~{Guti\'errez} $^{1,}$*\orcidA{}, Daniel~{G\'omez} $^{1,2}$\orcidC{}, Javier~{Castro} $^{1,2}$\orcidD{} and Rosa~Espínola~$^{1,2}$\orcidE{}}

\AuthorNames{Inmaculada  Guti\'errez, Daniel G\'omez,  Javier  Castro and Rosa Espínola}
\AuthorCitation{G\'omez, D.; Castro, J.; Guti\'errez, I.  and Esp\'inola, R.}

\AuthorCitation{Guti\'errez, I.; G\'omez, D.; Castro, J.; Esp\'inola, R.}

\address{%
$^{1}$ \quad Facultad de Estudios Estadísticos, Universidad Complutense de Madrid, Madrid, 28040, Spain 
, dagomez@estad.ucm.es (D.G.); jcastroc@estad.ucm.es (J.C.); rosaev@estad.ucm.es (R.E.) 
\\ 
$^{2}$ \quad Instituto de Evaluación Sanitaria, Universidad Complutense de Madrid, Madrid, 28040, Spain}

\corres{Correspondence: inmaguti@ucm.es}

\abstract{On the basis of network analysis, and within the context of modeling imprecision or vague information with fuzzy sets, we propose an innovative way to analyze, aggregate and apply this uncertain knowledge into community detection of real-life problems.  This work is set on the existence of one (or multiple) soft information sources, independent of the network considered, assuming this extra knowledge is modeled by a vector of fuzzy sets (or a family of vectors). This information may represent, for example, how much some people agree with a specific law, or their position against several politicians. We emphasize the importance of being able to manage the vagueness which usually appears in real life because of the common use of linguistic terms. Then, we propose a constructive method to build fuzzy measures from fuzzy sets. These measures are the basis of a new representation model which combines the information of a network with that of fuzzy sets, specifically when it comes to linguistic terms. We propose a specific application of that model in terms of finding communities in a network with additional soft information. To do so, we  propose an efficient algorithm and measure its performance by means of a benchmarking process, obtaining high-quality results. }



\begin{document}


\section{Introduction and Related Work}\label{sec:intro}

Social Network Analysis (SNA) is described as the study and  understanding of the relationships between two or more items. As one of the hottest topics of SNA, the Community Detection Problem (CDP) has become a problem of great interest  in modern statistics with applications in several fields \cite{fortunato,biological,basedOverlap}. 


Most of the algorithms and definitions of community detection problems assume that the only information available for identifying the clusters/communities in a network is the graph {which describes its structure}. This graph  can be non-directed and binary (all the relations in the network are equals) in the classical and most studied community detection problem \cite{biological}; non-directed and valued (the network is modeled by a weighted non-directed graph), or the case in which the relations are not symmetrical \cite{fortunato,directedMod,improved}. There are other interesting approaches focused on the incorporation of additional information to the crisp graphs, specifically to find communities in a network \cite{fumanalAffinity,fumanalWars}. Nevertheless, any of these approaches only considers the community detection problem  from a topological point of view, with a focus on the problem from the relations between nodes, but not considering other types of information that could be relevant in order to find communities in a real problem.

We illustrate our idea with an example. Let us present a situation in which we have a set $V$ of nodes which represents the members of a parliament, whose friendly relations are known to us by  the crisp graph $G=(V,E)$. Let us assume that the reason why they are interacting is because they are voting on a specific law in parliament.    This information (the voting problem)  and also their political preference on the law (or their capacity in the voting problem) could be relevant information  to identify the clusters in the network. 


To deal with this type of problem, in \cite{escim19,wcci2020,ipmu2020,infus,multipleBip}, the authors introduce a new element to the community detection problem: a capacity measure that tries to model and reflect the reason why the nodes are interacting in the network in addition to the interests of the nodes to remain united. From this perspective, in \cite{infus}, we present an efficient algorithm for a community detection problem that deals with networks and fuzzy measures in that sense. Furthermore, in \cite{ampliinfus}, we present a constructive method to build a $1$-additive fuzzy measure from a crisp valuation of the nodes in the network.


\textls[-15]{Nevertheless and due to the natural uncertainty in real problems, the information associated with the network nodes is not usually assumed to be crisp in a natural way. Uncertainty is associated with the lack of knowledge about the occurrence of some event. Within the last decades,  two important models are proposed to represent different types of uncertainty: randomness and vagueness/imprecision. Whereas the randomness emerges due to the lack of knowledge about the occurrence of some event \cite{novak2012mathematical},  vagueness} a phenomenon rises when trying to group together objects that share a certain property. A typical vague property is ``to be a small number'' or ``to be a tall person'', or (taking the previous example of the voting system in a parliament) ``to be against a specific law''.  In this way, the fuzzy linguistic approach has been successfully applied  to many problems \cite{herrera20002}.
    Taking into account this type of information, an important goal of this work is to provide a methodology to face community detection problems in networks with additional soft information. With the aim to extend some of the definitions and algorithms presented in \cite{ampliinfus} for crisp information, in this paper, we work on the basis of the existence of a vector (or a family of them) whose elements are no longer crisp values, but they are fuzzy sets that provide some type of soft information related to the individuals in the problem. In this context, another important objective of this work arises: we characterize a new representation tool which generalizes other existing models in the literature, regarding the nature of the information: the extended fuzzy graph based on a fuzzy vector (EFVFG). It is defined on the basis of a crisp networks and a vector of fuzzy sets. Another goal is to  extend it to a more complex scenario in which there is not only one type of information but many; in this situation, we strongly recommend consideration of the multi-dimensional extended fuzzy graph fuzzy vector-based (MEFVFG), which is defined on the basis of a family of vectors of fuzzy sets.

Then, we suggest a specific application of the new representation model, which is useful to obtain realistic partitions in a network with additional soft information. We present a competitive algorithm which introduces fuzzy sets to the process of grouping individuals. It is a modification of the well-known Louvain algorithm for crisp networks \cite{blondel} that allows us to deal with soft information in the network, which is developed on the basis of MEFVFG.  To guarantee the quality of the proposed methodology, we dedicate an important part of this work to its evaluation. The computational results showed in this work,  obtained through a benchmarking process  developed on the basis of some trapezoidal fuzzy sets, allow us to assert the good performance of our algorithm.

This paper is organized as follows. {In Section \ref{sec:prel}, we lay the foundations of the work, showing several concepts and definitions that are useful for the understanding and follow-up of the work. In Section \ref{sec:modelDef}, we characterize a new model representation based on soft information about the individuals of a network given by several fuzzy sets. After that, in Section \ref{sec:application}, we propose a specific application of that new tool, related to the community detection problem with additional soft information, which is a very live issue in the field of SNA. In order to evaluate the performance of the proposed methodology, we show some computational results in Section \ref{sec:computational}. The paper ends in Section \ref{sec:conclu} with some conclusions and a final discussion.}

\section{Preliminaries}\label{sec:prel}

\subsection{Fuzzy Sets}

Fuzzy sets were introduced by Zadeh as an extension of usual concept of set, and they have been applied in several fields \cite{bus,martinezchinos,Timothy}.

\begin{Definition}[Fuzzy set \cite{zadeh1965}]
Let $X$ denote a set. A fuzzy set in $X$, denoted by $\widetilde{A}$, is a set uniquely characterized by its membership function, defined by $\eta_A:X \longrightarrow [0,1]$ where, for every point $x\in X$,  $\eta_A(x)$ defines $x$'s ``grade of membership''.
\end{Definition}

In this work, we will focus on fuzzy sets over positive real numbers, so from now on, we will assume that $X=\mathbb{R^+}$ .

\color{black}

\begin{figure}[H]
\begin{minipage}{0.33\textwidth}
    \includegraphics[scale=1.5]{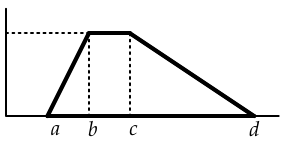}
    \end{minipage}
\caption{Graph 
 $G=(V,E)$ and fuzzy linguistic variable $\widetilde{L}$.}\label{cap3:fig:trapezoidal}
\end{figure}

Introduced by Zadeh \cite{zadeh1975} and  applied to the resolution of many real problems, the fuzzy linguistic variables were defined in situations in which imprecision or vagueness of a quantitative variable  are given in linguistic terms.  For example, a linguistic variable $\widetilde{L} = \{ L_1, \ldots, L_k\}$ can be characterized by $k$  membership functions, that is, the collection of its linguistic values, $U \subset \mathbb{R^+}$ is a universe of discourse and the meaning of each linguistic value is characterized by $\eta_{L_i}: U \longrightarrow  [0,1]$, which associates  each $u \in U$ with its compatibility. In the computational results section of this work, we consider a specific type of fuzzy sets that are commonly used to model the linguistic terms of a fuzzy linguistic variable: the trapezoidal fuzzy sets, whose shape is similar to Figure \ref{cap3:fig:trapezoidal}.

	\begin{Definition}[Trapezoidal fuzzy set \cite{rankingTriangular}]\label{cap3:def:TrapfuzzyNumber} The fuzzy set ${A}=\left(a,b,c,d\right)$ is said to be trapezoidal if its membership function$\eta_{A}$ is defined by:
%

\begin{equation}\eta_{\mathcal{A}}
\left(x\right)=\begin{cases} 0 &  \textup{ if }  x<a \\
\frac{x-a}{b-a} & \textup{ if } a\leq x\leq b \\
1 & \textup{ if }  b< x\leq c \\
\frac{x-d}{c-d} & \textup{ if }  c< x\leq d \\
0 & \textup{ if } d<x 
\end{cases}\end{equation}
	\end{Definition}

\vspace{-12pt}

In this work, some of the input values are defined as fuzzy sets. Nevertheless, to face the final goal, we need to ``convert'' that soft information into crisp values. This process of obtaining a single output from the output of an aggregated fuzzy set is known as ``defuzzification'', and there are different methods to carry on with it \cite{defuzzCriteria,defuzz:delphi}. In general, a formal definition of a defuzzification operator is shown below.

\begin{Definition}[Defuzzification operator \cite{defuzz:delphi}]
Given a universe $X$, the operator $D:F(X)\rightarrow X$ which maps the fuzzy sets on $X$ into elements of the $X$ is said to be a defuzzificatin operator.
\end{Definition}

\subsection{Networks with Additional Information: An Algorithm That Deals with CDP and Fuzzy~Measures}\label{sec:CDPFuzzy}

A graph is a pair $G=(V,E)$, in which $V=\{1,2,\dots , n\}$ is a set of individuals called nodes or vertices, and $E=\{\{i,j\}| i,j\in V\}$ is an unordered set of pairs of nodes called edges or arcs. A graph is unequivocally defined by its adjacency matrix $A$, which is characterized as $\left(A\right)_{i,j\in V}=1$ if $\{i,j\}\in E$, and $0$ otherwise. A graph is said to be valued or weighted if there is a function $w:E\rightarrow \mathbb{R}$ which assigns a weight to each edge. In this type of graph, the adjacency matrix not only represents the existence of an edge between two nodes but also shows the weight of each edge by displaying the value assigned to each edge by the $w$ function.

The community detection problem is an important problem in the SNA field, the goal of which is to find a ``good'' partition of the set of nodes. A partition is considered to be only as good as how internally homogeneous and externally heterogeneous the defined groups are in terms of the connections between individuals. The modularity measure defined in \cite{biological}, usually denoted by $Q$, is a quality function of the partitions which somehow measures the strength of the division of a graph in a partition of communities. $Q$  is usually considered as a function to be maximized.

\begin{Definition}[{Modularity} \cite{girvanNewman}]\label{cap2:def:modularity}
Let $G=\left(V,E\right)$ denote a graph with adjacency matrix $A$. Let  $i,j\in V$ and $m=|E|$.  The modularity function of the partition of $V$, $P$, is characterized by

	\begin{equation}\label{eq:modularity}
		Q\left(P\right)= \frac{1}{2m} \sum_{i,j\in V} \left[{A}(i,j)- \frac{k_i \ k_j}{2m}
	\right] \delta(C_i, \ C_j)
	\end{equation} 
	in which $k_i$ is the degree of $i$ and  $C_i$ is the group to which $i$ is assigned; $\delta\left({C_i, C_j}\right)=1$  if $C_i=C_j$, and $\delta\left({C_i, C_j}\right)=0$ otherwise.
\end{Definition}

Without detriment regarding the worthiness of classic approaches, some authors agreed on the importance of including as much information as possible in the network analysis process, regardless of the direct crisp connections between individuals defined by the edges. We find several approaches with a common idea: the more information is considered, the more realistic the results obtained, either in terms of partitions or any other notion \cite{centralityPower,fumanalAffinity,edgeBet}. Specifically, this work is set on the basis of the idea introduced in \cite{infus}. In that preliminary work, the authors proposed a methodology to find realistic  communities in a graph in terms of a fuzzy measure, defining some additional information about the synergies between the individuals. That method was based on the Louvain algorithm \cite{blondel} with a main difference: the calculation of modularity not only considers the adjacency matrix but also some additional information matrix, specifically, one obtained from the mentioned affinity fuzzy measure. This methodology, named \textit{Duo Louvain}, is summarized in Algorithm \ref{alg:duoLouvain}. The main difference with respect to the common Louvain method can be seen in line $15$ of the pseudo-code: the variation of modularity obtained when moving the node $o^i$ to the community to which its neighbor $e_j$ belongs, $\Delta Q_{o^i}(e_j)$, is calculated in any matrix $M$, which is different from the adjacency matrix. This methodology was adapted to different scenarios in later works \cite{multipleBip, polarizacion:covid,wcci2020} related to a variety of fuzzy measures. Let us emphasize that this methodology is far more powerful than being limited to the consideration of fuzzy measures. It can also be considered in any  other scenario beyond them, provided that any additional information that can be aggregated, in any form, into a matrix, is available.

As mentioned, a quick overview of that methodology is showed in Algorithm \ref{alg:duoLouvain}, where $\pi(V)$ denotes all the feasible permutations of the elements of $V$; $o=(o^1,\dots ,o^n)\in \pi(V)$ is one of these orders; $H(o^i)$ denotes the set of neighbors of $o^i\in V$ or, what is the same, the nodes with which $o^i$ is directly connected, and $\Delta Q_{o^i}(e_j)$ denotes the variation of modularity obtained when moving $o^i$ to the community to which $e_j$ belongs.
\begin{algorithm}[H]
\begin{algorithmic}[1]
 \STATE\textbf{Input}: $\left(A, \ {M}\right)$;
\STATE\textbf{Output}:
$P$;
\STATE\textbf{Preliminary}
\STATE $C_i \leftarrow \{i\},$ $\forall i \in V$ \ \  (each node $i$ is an isolated community);
\STATE $P \leftarrow \left(1, 2,\ldots, n\right)$ \ \ (initial partition);
\STATE\textbf{end Preliminary}
 \STATE\textbf{Phase 1}
 \STATE Take $o=\left(o^1, \dots, o^i,\ldots, o^n\right) \in \pi(V)$;
 \STATE $stop \leftarrow 0$;
\WHILE{$(stop==0)$}
\STATE{$stop \leftarrow 1$} 
\FOR{$(i=1)$ \TO $(n)$}
 \STATE $\left(e_1,\dots, e_h\right) \leftarrow H(o^i)$ \ \ (find the neighbors of $o^i$ in $A$);
\FOR{$(j=1)$ \TO $(h)$}
\STATE Calculate $\Delta Q_{o^i}(e_j)$ in ${M}$;
\ENDFOR
 \vspace{0.8mm}
\STATE ${\displaystyle j^* \leftarrow \Big\{e_{\ell} \ | \ \Delta Q_{o^i}(j^*)=\max_{\ell\in\{1 \dots , h\}}\big\{\Delta Q_{o^i}(e_{\ell})\big\}\Big\}}$; 
\IF{$(\Delta Q_{o^i}(j^*)>0)$}
\STATE $C_{P\left(o^i\right)} \leftarrow C_{P\left(o^i\right)}\backslash{\{o^i\}}$; 
\vspace{0.8mm}
\STATE $C_{P\left(j^*\right)} \leftarrow C_{P\left(j^*\right)} \cup \{o^i\}$; 
\vspace{0.8mm}
\STATE $P\left(o^i\right) \leftarrow P\left(j^*\right)$;
\vspace{0.8mm}
\STATE $stop \leftarrow 0$;
\ENDIF
\ENDFOR
\ENDWHILE
\STATE\textbf{end Phase 1}
 \STATE\textbf{Phase 2} 
\STATE Calculate $A^*$ from  $A$ (nodes of $A^*$ are the communities previously found in $A$);
\STATE Calculate $M^*$  from  $M$ (nodes of $M^*$ are the communities previously found in $M$);
\IF {$(A^* \neq A)$}
\STATE $A \leftarrow A^*$;
\STATE $M \leftarrow M^*$;
\STATE Apply \textbf{Phase 1} and \textbf{Phase 2};
\ENDIF
\STATE\textbf{end Phase 2}
\STATE\textbf{return}$\left(P\right)$;

\end{algorithmic}\caption{\textit{Duo Louvain}}
\label{alg:duoLouvain}
\end{algorithm}

\color{black}

\vspace{-12pt}

\section{Model Definition: Building Extended Fuzzy Graphs from Graphs with Fuzzy Nodes~Information}\label{sec:modelDef}

In this section, we work on the definition of a new representation tool. Firstly, we do this in a uni-dimensional scenario, assuming there is an additional {fuzzy} information vector related to the individuals of a set $V$, denoted by $\widetilde{f}=\left(\widetilde{f_1}, \dots ,\widetilde{f_r}\right)$.  For each $i \in V$, the fuzzy set $\widetilde{f_i} $ (characterized by its membership function $\eta_{f_i}$) represents the vague or imprecise information associated to the node $i$ of some characteristic or evidence. This fuzzy modelization is especially useful (but not only) when the information associated with each node is gathered (for example) by a linguistic term.  Specifically, in this case, we could work with linguistic terms $\widetilde{f_i} \in \widetilde{L} $. By analogy with \cite{wcci2020}, we first propose a characterization of a fuzzy Sugeno $\lambda$-measure from this fuzzy vector $\widetilde{f}$. This measure is denoted by $\mu_{f,p}$.

\begin{Definition}[{Fuzzy Sugeno $\lambda$-measure obtained from fuzzy sets}]\label{def:mufp}
Given the set \linebreak \mbox{$V=\{1,2,\dots ,n\}$}, let $\widetilde{f}=\left(\widetilde{f_1}, \dots ,\widetilde{f_n}\right)$ denote a vector of fuzzy sets defined over a universe $U \subset \mathbb{R}^+ $  (i.e., $\eta_{f_i}: U \longrightarrow [0,1]$), and let $D:F(\mathbb{R}^+)\rightarrow \mathbb{R}^+$ denote a defuzzification operator. Then, for any $p \in (0,1]$ and $i\in V$, a natural definition of is $\mu_{f,p}$ is:

\begin{equation}\label{eq:mufp}\mu_{f,p}(i)=\frac{pD( \widetilde{f_i})}{\sum_{k=1}^n D(\widetilde{f_k})}, \hspace{5mm} \forall i \in V  \end{equation}
 \textup{ where } $\mu_{f,p}(M\cup N)=\mu_{f,p}(M)+\mu_{f,p}(N)+\lambda\mu_{f,p}(M)\mu_{f,p}(N),
\ \forall M, N \subseteq V$, \textup{with} $M\cap N=\emptyset  \textup{ and }\lambda+1=\prod\limits _{i=1}^n(1+\lambda\mu_{f,p}(i)), \textup{ being  } p\in(0,1]$. 
\end{Definition}

Note that the interpretation of $\mu_{f,p}$ depends on $p$. Specifically, 

\begin{Proposition}\label{prop:defAdd}
Given the parameter $p=1$, the function $\mu_{f,p}$ is a fuzzy Sugeno $\lambda$-measure $1$-additive.
\end{Proposition}


\begin{proof}
Because of the properties of the Sugeno $\lambda$-measures \cite{phdSugeno} in addition to the assumption of $p=1$, we have $\lambda=0$. Then, $\forall M\subseteq V$, $\mu_{f,p}(M)=\frac{\sum_{\ell\in M} D(\widetilde{f_{\ell}})}{\sum_{k=1}^n D(\widetilde{f_k})}$, so $\mu_{f,p}$ meets the conditions of fuzzy measures \cite{sugeno}, Sugeno $\lambda$-measures \cite{phdSugeno} and $1$-additivity \cite{gra2}.

\begin{itemize}
    \item $\mu_{f,p}(\emptyset)=0$ Trivial.
    \item $\mu_{f,p}(V)=\frac{D(\widetilde{f_1})}{\sum_{k=1}^n D(\widetilde{f_k})} + \dots + \frac{D(\widetilde{f_n})}{\sum_{k=1}^n D(\widetilde{f_k})}=1$.
    \vspace{4mm}
    \item Let $M\subseteq N\subseteq V$.
    Then, $\mu_{f,p}(N)=\frac{\sum_{\ell \in N}D(\widetilde{f_{\ell}})}{\sum_{k=1}^n D(\widetilde{f_k})}=$
    $\frac{\sum_{\ell \in M}D(\widetilde{f_{\ell}}) +     \sum_{t \in N\backslash M}D(\widetilde{f_t})}{\sum_{k=1}^n D(\widetilde{f_k})} \geq \frac{\sum_{\ell \in M}D(\widetilde{f_{\ell}})}{\sum_{k=1}^n D(\widetilde{f_k})} = \mu_{f,p} (M)$,  so $\mu_{f,p}$ is a fuzzy measure.
    \item \textbf{Sugeno $\lambda$-measure}. Trivial by definition.
    \item \textbf{$1$-additivity}: regarding \cite{gra2}, it is trivial if $\forall i \in \{1, \dots , n\}$, we define $a_i=\mu_{f,p}(i)$.
    \end{itemize}
\end{proof}

\begin{Proposition}\label{prop:mufpGeneral}
Given the parameter $p\in(0,1)$,  $\mu_{f,p}$ is a fuzzy Sugeno $\lambda$-measure.
\end{Proposition}
\begin{proof}
The proof is similar to that of Proposition \ref{prop:defAdd}.
\end{proof}

So, we generalize the notion of {\it extended fuzzy graph vector based}, $\widetilde{G}=\left(V,E,\mu_{x,p}\right)$ \cite{wcci2020} to a scenario where the additional information is not provided by a crisp vector $x$, but it comes from a vector of fuzzy sets, $\widetilde{f}=\left(\widetilde{f_1}, \dots ,\widetilde{f_r}\right)$.

\begin{Definition}[{\bf Extended fuzzy graph fuzzy vector based (EFVFG)}]
Let $G=(V,E)$ denote a graph with $n=|V|$ individuals and $m=|E|$ edges. Let  $\widetilde{f}=\left(\widetilde{f}_1,\dots,\widetilde{f}_n\right)$ denote a vector of {fuzzy sets in membership function form, each of them related to an individual of $V$.}  Let $D:F(\mathbb{R^+})\rightarrow \mathbb{R^+}$ denote a defuzzification operator, and given the parameter $p\in(0,1]$, let $\mu_{f,p}$ denote the fuzzy Sugeno $\lambda$-measure obtained from $\widetilde{f}$. Then, the tuple $\widehat{G}=\left(V,E,\mu_{f,p}\right)$  is said to be a fuzzy extended graph based on the fuzzy vector $\widetilde{f}$.
\end{Definition}

\begin{Example}\label{ex1} Let $\widetilde{L}=\{Very-Low,Low, Medium, High, Very-High\}$ denote a fuzzy linguistic variable defined over the universe $U=[0,100]$, which is  characterized by the corresponding membership functions {$\eta_{VL}, \eta_{L}, \eta_{M}, \eta_{H}, \eta_{VH}: [0,100] \rightarrow [0,1]$} \textls[-15]{associated with the different linguistic terms that represent how in agreement a person is with some law denoted by} $LW^1$. Let $G=(V,E)$ define a {\it cyclic} graph with  $V=\{1,2,3,4,5,6,7,8\}$ and  $E= \{ (1,2), (2,3), (3,4), \ (4,5), (5,6), \\ (6,7), (7,8), (8,1) \}$, and finally, let  $\widetilde{f}= \left(  \widetilde{f_1}, \ldots, \widetilde{f_8} \right)=(VL,VL,L,VL,H, VH, H, VH)$ denote a vector of fuzzy sets that models the linguistic terms affinity of these $eight$ nodes of $V$ to the law $LW^1$.

From the previous definition, it is possible to build  (for any $p \in (0,1]$) the extended fuzzy graph associated with the fuzzy vector $\widetilde{f}$ and the graph $G=(V,E)$, that is: $\widehat{G}=\left(V,E,\mu_{f,p}\right)$.

\end{Example}

\vspace{-12pt}

Assuming we can have more than one characteristic associated with each node in a network, we go beyond the uni-dimensional case by considering there is not only a {vector of fuzzy sets} $\widetilde{f}$, but a family of them, $\left( \widetilde{f^1}, \dots ,\widetilde{f^r}\right)$, each of them defining some extra knowledge about the individuals.

\color{black}

\begin{Definition}[{\bf Multi-dimensional extended fuzzy graph fuzzy vector based (MEFVFG)}]
Let $G=(V,E)$ denote a graph with $n$ nodes and $m$ edges, and let  $\left(\widetilde{f}^1,\dots,\widetilde{f}^r\right)$ denote a family of $r$ independent vector of $n$ fuzzy sets,  each of them defining a type of information, so that $\forall \ell=1,\dots ,r; i=1,\dots,n$ the component $\widetilde{f^{\ell}_i}$ is the fuzzy set {related to the characteristic $\ell$ and the individual $i$} with membership function $ \eta_{f^{\ell}_i}: U^{\ell} \longrightarrow [0,1]$.  Let $D:F(\mathbb{R^+})\rightarrow \mathbb{R^+}$ be a defuzzification operator (that we will assume that is the same for all characteristics $\widetilde{f}^\ell$), and  let $p^{\ell}\in(0,1]$ be a parameter in $(0,1]$.

Then, the tuple $\widehat{G}=\left(V,E,\left(\mu_{f^{1},p^{1}},\dots , \mu_{f^{r},p^{r}}\right)\right)$  is said to be a multi-dimensional extended fuzzy graph (MEFVFG) based on the $r$ fuzzy vectors ($\widetilde{f^1}, \ldots, \widetilde{f^r}$ ). 

\end{Definition}

\begin{figure}[H]
\begin{minipage}{0.23\textwidth}
    \includegraphics[scale=0.5]{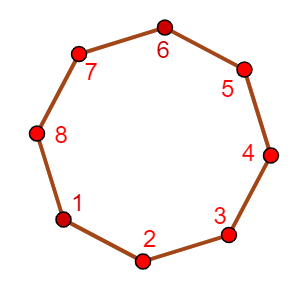}
    \end{minipage}\begin{minipage}{0.7\textwidth}
    \includegraphics[scale=0.5]{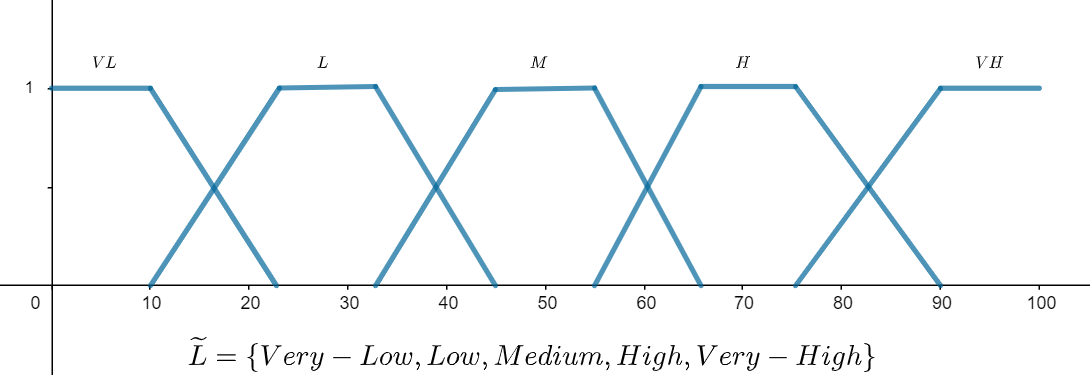}
\end{minipage}\caption{Graph $G=(V,E)$ and fuzzy linguistic variable $\widetilde{L}$.}\label{fig:ex}
\end{figure}

\begin{Example} Let $\widetilde{L}=\{Very-Low,Low, Medium, High, Very-High\}$ denote a fuzzy linguistic variable defined over the universe $U=[0,100]$ characterized by the  corresponding membership functions {$\eta_{VL}, \eta_{L}, \eta_{M}, \eta_{H}, \eta_{VH}: [0,100] \rightarrow [0,1]$} associated with the different linguistic terms that represent how in agreement a person is with a specific law, $LW^1$. Let $G=(V,E)$ denote a {\it cyclic} graph with $V=\{1,2,3,4,5,6,7,8\}$ and $E= \{ (1,2), (2,3), (3,4),(4,5), (5,6), (6,7), (7,8), \\ (8,1) \}$, and  finally let $\widetilde{f^1}= \left(  \widetilde{f^1_1}, \ldots, \widetilde{f^1_8} \right)=(VL,VL,L,VL,H, VH, H, VH)$ be a vector of fuzzy sets that models the linguistic terms affinity of these eight nodes to the congress proposal $LW^1$ and let $\widetilde{f^2}= \left(  \widetilde{f^2_1}, \ldots, \widetilde{f^2_8} \right)=(M,M,M,M,L, L, L, H)$ be a vector of eight fuzzy sets that models the linguistic terms affinity of the eight nodes to the congressional bill, which is denoted by $LW^2$ (see Figure \ref{fig:ex}). 

From the previous definition, for any $p_1$,$p_2\in (0,1]$,  it is possible to build the {MEFVFG} associated with the two fuzzy vectors $(\widetilde{f^1}, \widetilde{f^2})$ as

$$\widehat{G}=\left(V,E,\left(\mu_{f^1,p_1}, \mu_{f^2,p_2} \right)\right).$$ 

\end{Example}

Let us remark that this last case generalizes other existent tools, such as for example fuzzy graphs defined in \cite{rosenfeldFG}, which actually only define relations between connected individuals, the extended fuzzy graphs \cite{infus}, in which the additional information is about the relations between the elements but not on the individuals itself; or the  \cite{wcci2020} whose additional information is about individuals, but it is crisp.

\section{An Application: Social Network Analysis with Soft Information}\label{sec:application}

As a specific application of the new proposed model, we take up  firstly the idea  introduced in \cite{wcci2020} about community detection in graphs taking into account the information given by a crisp vector, $x$. In that preliminary work, we set up the philosophy of finding groups in a network when there is a vector of crisp values providing some additional information. Now, we generalize this idea, starting from extended fuzzy graphs that are built from $r$ fuzzy vectors $\widetilde{f^1}, \ldots, \widetilde{f^r}$.

To find ``good'' communities in $\widehat{G}$, we have to extend the  \textit{Sugeno Louvain} Algorithm described in \cite{wcci2020} to a multi-dimensional stage with more than one vector of additional information, with the peculiarity that the components of the vectors considered are no longer crisp values but fuzzy sets: therefore, actually, what we have is an  MEFVFG. We illustrate the problem of community detection based on an MEFVFG in Example \ref{ex:comMEFVFG}.

\vspace{-2pt}

\begin{Example}\label{ex:comMEFVFG}
\textls[-25]{We consider a chain with $12$ nodes, represented by the crisp graph $G=(V,E)$ (see Figure \ref{fig:my_label}) about which we have additional information, $(\widetilde{f}^1, \widetilde{f}^2,\widetilde{f}^3, \widetilde{f}^4)$, and the defuzzification operator $D$}, so $(D(\widetilde{f}^1)=(\mathbf{9,  9.5, 10}, 1, 0.5, 1, \mathbf{9.5, 8, 10},  1, 1, 2), \ D(\widetilde{f}^2)=(\mathbf{10, 9.5, 9}, 1, 0.5, 1, \mathbf{9, 9, 9.5}, 1.5, 2, 0.5), \\ D(\widetilde{f}^3)=(\mathbf{9.5,8.5, 10}, 1.5, 1, 1, \mathbf{10, 9, 9.5}, 0.9, 1, 1),  \ D(\widetilde{f}^4)=(\mathbf{9, 9, 10}, 1, 1, 1, \mathbf{10, 9.5, 9}, 0.5, 1, 1))$. \textls[-25]{These fuzzy sets represent the opinion of $12$ people about $four$ different films. We accept that there are more synergies between those people who have similar preferences.  Partition} $P=\{\{1,2,3,4\},\{5,6,7,8\}, \{9,10, \\ 11,12\}\}$ is obtained with any     algorithm based on modularity optimization.  Nevertheless, if the additional information is considered, the partition provided by the {\textit{Multi-Dimensional fuzzy Sugeno--Louvain $1$-additive}} is $P^f=\{\{1,2,3\}, \{4,5,6\},\{7,8,9\}, \{10,11,12\}\}$.

\vspace{-9pt}
\begin{figure}[H]
    \includegraphics[scale=0.6]{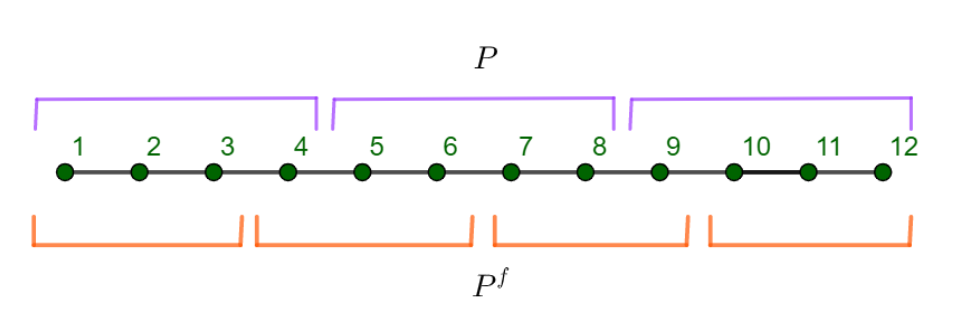}
    \caption{Chain with 
 $12$ nodes. Partitions $P$ and $P^f$.}
    \label{fig:my_label}
\end{figure}
\end{Example}

\vspace{-1pt}

\begin{figure}[H]
    \includegraphics[scale=0.7]{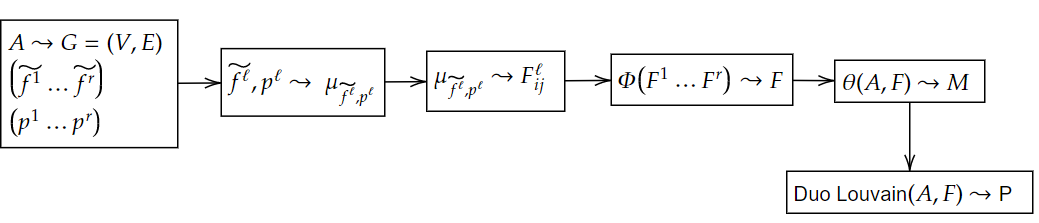}
\caption{Flowchart of the methodology 
 \textit{Multi-dimensional Fuzzy Sugeno--Louvain}.}\label{fig:flow}
\end{figure}

The proposed method, named \textit{ Multi-dimensional Fuzzy Sugeno Louvain}, is based on the Louvain algorithm \cite{blondel}. The main point is to summarize all the knowledge of the  MEFVFG into two matrices: $A$ that represents the direct connections between the nodes (edges), and $F$ summarizes the additional information given by the family of vectors of fuzzy sets $\left( \widetilde{f^1}, \dots , \widetilde{f^r}\right)$. The weighted graph associated to a fuzzy Sugeno $\lambda$-measure \cite{wcci2020} $\mu_{f^\ell,p^{\ell}}$ is essential, and it is considered in terms of a multi-dimensional scenario (MAWG) (one weighted graph with adjacency matrix $F^{\ell}$ associated to each $\mu_{f^{\ell},p^{\ell}}$). This methodology to find realistic partitions in an MEFVFG explained below is summarized in Algorithm \ref{cap3:alg:SugenoLouvainMulti}, which includes its pseudo-code, and in Figure \ref{fig:flow}, which shows a flowchart of the process.

\color{black}
\begin{itemize}
 \item \textbf{\textls[-25]{Step 1: definition of the MAWG}.} \textls[-25]{Given the fuzzy Sugeno} $\lambda$-measures $\left(\mu_{f^1,p^1}, \dots , \mu_{f^r,p^r}\right)$ obtained from  $\left(\widetilde{f}^1, \dots \widetilde{f}^r\right)$ and $\left(p^1,\dots ,p^r\right)$, and the defuzzification operator $D$, \linebreak matrices  $\left({F}^{1}, \dots, {F}^{r}\right)$ are calculated as
 \begin{equation}\label{eq:calculoF}
    {\small F_{ij}^{\ell}=\phi \left(Sh_i(\mu_{f^{\ell},p^{\ell}})-Sh_i^j(\mu_{f^{\ell},p^{\ell}}), Sh_j(\mu_{f^{\ell},p^{\ell}})-Sh_j^i(\mu_{f^{\ell},p^{\ell}})\right) } 
\end{equation}
being $\phi:[-1,1] \rightarrow [0,1]$ a bi-variate aggregation operator \cite{orderedWeights}; $Sh_i(\mu_{f^{\ell},p^{\ell}})$ and $Sh_i^j(\mu_{f^{\ell},p^{\ell}})$ the Shapley values of $i$ on $\mu_{f^{\ell},p^{\ell}}$ in the presence of all the elements of $V$ or $V\backslash\{j\}$, respectively \cite{shapley}.
 
 \item \textbf{Step 2:  information aggregation.} Matrices  ${F}^{1}, \dots ,{F}^{r}$ are aggregated to obtain  the matrix $F$. The aggregation function  $\Phi:\Pi^r \rightarrow \Pi$ is used, being $\Pi$ the set of quadratic $n$-matrices. Particularly, we suggest the use of a matrix aggregator based on the classical aggregation operators with element to element transformation:
 $ F=\Phi\left({F}^1, \dots, {F}^{r}  \right)$. 
\end{itemize}

After this aggregation process, the method \textit{Duo Louvain} has to be applied \cite{ampliinfus,multipleBip}, considering the matrix $M=\theta\left(A,F\right)$, being $\theta:\Pi^2\rightarrow \Pi$ an aggregation function. That method can consider the information of two matrices when finding communities in a graph.

\begin{algorithm}[H]
 \begin{algorithmic}[1]
\STATE\textbf{Input}: $\left(A, \left(\ \widetilde{f}^1, \dots , \widetilde{f}^r\right),\left(\ p^1, \dots , p^r\right)\right)$, $A$ represents $G=(V,E)$; $\widetilde{f}^{\ell}$ is a vector of fuzzy sets; $p^{\ell}\in [0,1)$, $\forall \ell = 1, \dots , r$;
\STATE\textbf{Output}: $P$;
\STATE\textbf{Preliminary}
\FOR{$(\ell=1)$ \TO $(r)$}
\STATE Calculate $\mu_{f^\ell,p^\ell}$ (fuzzy Sugeno $\lambda$-measure from $\widetilde{f}^{\ell}$);
\STATE  ${F}_{ij}^\ell \leftarrow \phi \left(Sh_i(\mu_{f^\ell,p^\ell})-Sh_i^j(\mu_{f^\ell,p^\ell}),  Sh_j(\mu_{f^\ell,p^\ell})-Sh_j^i(\mu_{f^\ell,p^\ell})\right)$, $ $ $\forall i,j\in V$;
\ENDFOR
\STATE $F \leftarrow \Phi\left({F}^1, \dots , {F}^r\right)$;
\STATE $M \leftarrow\theta\left(A,F\right)$;
\STATE\textbf{end Preliminary}
\STATE $P \leftarrow$ \textit{Duo Louvain}$\left(A,M\right)$;
\STATE\textbf{return}$\left(P\right)$;
 \end{algorithmic}\caption{{\textit{Multi-dimensional Fuzzy Sugeno--Louvain}} }
 \label{cap3:alg:SugenoLouvainMulti}
\end{algorithm}

\vspace{-12pt}

\begin{Remark}\label{remark:agregaVector}
 The concept of ``what is a good group''  depends on the operator $\Phi$ applied.  In the case that $\Phi$ is a disjunctive operator, groups are composed by elements among which there are strong synergies regarding any evidence or characteristic (any fuzzy vector). The size of the groups that are somehow similar regarding the additional information will increase the more vectors are considered.
In contrast, where $\Phi$  is a conjunctive operator, the groups are composed by elements among which there are strong synergies in all the evidence or characteristics. The size of the groups that are somehow similar regarding the additional information will increase the less vectors are considered.
Particularly, we consider the most popular \textit{ordered weighted averaging aggregation operators}, OWA~\cite{owaYager}: maximum, minimum and average.
\end{Remark}


As in the uni-dimensional problem with crisp information addressed in \cite{wcci2020}, the exponential complexity concerning the calculation of the Shapley value may be avoided by considering an additive fuzzy measure. For this reason, in this paper, we suggest the specific characterization of $\mu_{f^\ell,p^\ell}$ when $p=1$. On this basis, as $\mu_{f^\ell}^a$ is a $1-$additive fuzzy measure \cite{gra2}, it holds:
$$Sh_i(\mu_{f^\ell}^a)=\frac{D(\widetilde{f}^{\ell}_i)}{\sum_{\begin{subarray}{1} k=1 \end{subarray}}^n D(\widetilde{f}^{\ell}_k)} \ \
\textup{ and } \ \ Sh_i^j(\mu_{f^{\ell}}^a)=\frac{D(\widetilde{f}^{\ell}_i)}{\sum_{\begin{subarray}{1} k=1 \\ k \not =j \end{subarray}}^n D(\widetilde{f}^{\ell}_k)}$$

In this context, we propose a specific application of the Algorithm {\textit{Multi-dimensional Fuzzy Sugeno--Louvain}}. For every $\ell=1,\dots, r$, the characterization of  $\mu_{f^\ell}^a$ only depends on the calculation of ${F}^\ell$. Then, the complexity of the method {\textit{1-additive Multi-dimensional Fuzzy Sugeno--Louvain}} is equal to that of the Louvain algorithm \cite{blondel}.

\begin{algorithm}[H]
 \begin{algorithmic}[1]
  \STATE\textbf{Input}: $\left(A, \ \left(\widetilde{f}^1, \dots , \widetilde{f}^r\right)\right)$, $A $ is a representation of $G=(V,E)$;  $\widetilde{f}^{\ell}$ is a vector of fuzzy sets, $\forall \ell = 1, \dots , r$;
\STATE\textbf{Output}: $P$;
\STATE\textbf{Preliminary}
\FOR{{$(\ell=1)$ \TO $(r)$}}
\STATE ${F}_{ij}^\ell \leftarrow \phi \{ |\frac{D(\widetilde{f}_i^\ell)}{\sum_{\begin{subarray}{1} k=1 \end{subarray}}^n D(\widetilde{f}_k^\ell)}-\frac{D(\widetilde{f}_i^\ell)}{\sum_{\begin{subarray}{1} k=1 \\ k \not =j  \end{subarray}}^n D(\widetilde{f}_k^\ell)}|, \ |\frac{D(\widetilde{f}_j^\ell)}{\sum_{\begin{subarray}{1} k=1 \end{subarray}}^n D(\widetilde{f}_k^\ell)}-\frac{D(\widetilde{f}_j^\ell)}{\sum_{\begin{subarray}{1} k=1 \\ k \not =i  \end{subarray}}^n D(\widetilde{f}_k^\ell)}|\}$;
\ENDFOR 
\STATE $F \leftarrow \Phi\left({F}^1, \dots , {F}^r\right)$;
\STATE $M \leftarrow \theta\left(A,F\right)$;
\STATE\textbf{end Preliminary}
\STATE $P \leftarrow$ \textit{Duo Louvain}$\left(A,M\right)$;
\STATE\textbf{return}$\left(P\right)$;
 \end{algorithmic}\caption{\textit{1-additive Multi-dimensional Sugeno--Louvain}}
 \label{cap3:alg:SugenoLouvainAditivaMulti}
\end{algorithm}

\begin{Example} We illustrate the idea of our methodology in a simple case considering there is  a vector of fuzzy sets. Let us consider the situation described in Example \ref{ex1} in which we have a cycle of eight nodes and the information associated to each node is described in linguistic terms 
$\widetilde{f}= \left(  \widetilde{f_1}, \ldots, \widetilde{f_8} \right)=(VL,VL,L,VL,H, VH, H, VH)$. The whole information is summarized in the Figure \ref{fig:ex2}. Now, let us assume that these  linguistic fuzzy variables defined over $U=[0,100]$ are modeled in terms of the following four fuzzy trapezoidal sets. 
$ \widetilde{VL}= (0,0,10,25) $, $ \widetilde{L}= (5,10,20,25) $, $\widetilde{M}= (30,40,60,70)$, $ \widetilde{H}= (60,70,80,100) $, $ \widetilde{VH}= (75,90,100,100) $.
It is possible to see that if we apply the {\textit{Fuzzy Sugeno--Louvain $1$-additive}}, just in the unidimensional so we only have one matrix $F^{\ell}$, to this extended fuzzy graph the partition obtained for any $p$ is $P^f=\{\{1,2,3,4\}, \{4,5,6,7,8\} \}$.
\begin{figure}[h]\centering
\begin{minipage}{0.23\textwidth}
    \includegraphics[scale=0.5]{Images/rueda.png}
    \end{minipage}\begin{minipage}{0.7\textwidth}
    \includegraphics[scale=0.5]{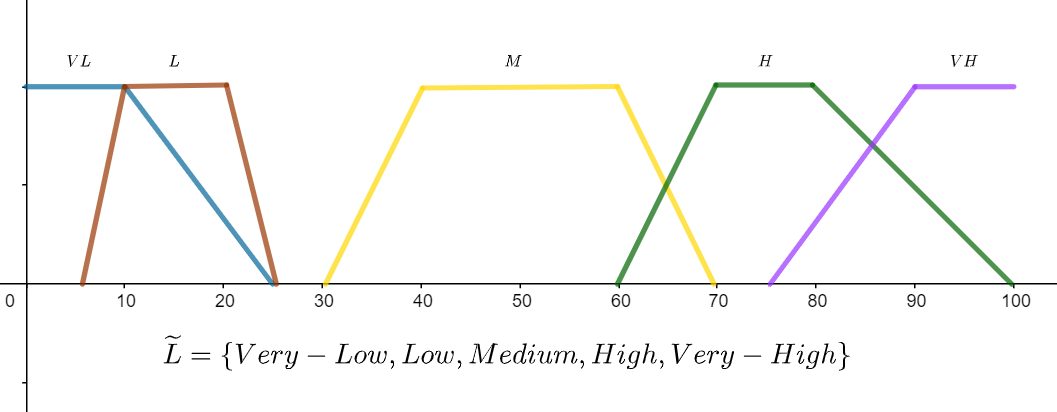}
\end{minipage}\caption{Graph 
 $G=(V,E)$ and fuzzy linguistic variable $\widetilde{L}$.}\label{fig:ex2}
\end{figure}

\color{black}

\end{Example}

\vspace{-17pt}

\section{Computational Results}\label{sec:computational}

When a new method is proposed, an evaluation of its performance is required. This process can be addressed comparing the results obtained with the method under evaluation with respect to other proposals established in the literature to solve the same problem. Nevertheless, in our case, as the community detection problem with additional soft information has never been faced before, we cannot compare our method with other proposals of the literature. Then, we work on an evaluation process. For this, we consider several reference models \cite{benchmarking2} to which we apply our methodology, whose performance is quantified with the calculation of the Normalized Mutual Information (NMI) \cite{nmi}.

\begin{Definition}[Normalized Mutual Information (NMI) \cite{nmi}]\label{cap2:def:nmi}
	\textls[-25]{Let $X=\{x_i\}_{i\in V}$ and} \mbox{$Y=\{y_i\}_{i\in V}$} denote two disjoint partitions of the graph $G=\left(V,E\right)$. 
	Let $P(x)$ denote the probability that a random node is assigned to the community $x$, and let $P(x,y)$ denote the conditional probability that a random node is assigned to the community $x$ in the partition $X$ and assigned to the community $y$ in the partition $Y$.
	The Shannon entropy of $X$ is calculated as $H(X)=-\sum_{x}P(x)\*log(P(x))$; and the Shannon entropy of $X$ and $Y$ is calculated as $H\left(X,Y\right)=-\sum_{x}\sum_{y}P(x,y)\*log(P(x,y))$. The  Mutual Information ($MI$) among partitions $X$ and $Y$ is defined as:
	\begin{equation}\label{cap2:eq:mutualInf}
	MI\left(X,Y\right)=\sum_{x}\sum_{y}P\left(x,y\right)\*log\frac{P\left(x,y\right)}{P\left(x\right)\*P\left(y\right)}
	\end{equation}
	
	Then, $NMI$  is a normalization of Equation \eqref{cap2:eq:mutualInf}.
	\begin{equation}\label{cap2:eq:nmi}
	NMI\left(X,Y\right)=\frac{2\*MI\left(X,Y\right)}{H(X)+H(Y)}
	\end{equation}
\end{Definition}

Although there can be some issues with the basic version of the measure \cite{AMELIO}, we consider this measure because, to the best of our knowledge, it is fair enough to compare how similar two partitions are; i.e., NMI allows us to quantify how much the partition provided by our method resembles the considered standard partition.

\color{black}

In all the benchmark models we present, there are two components: the adjacency matrix $A$ and the additional information matrix $F$, which are obtained from {some aggregation of a family of vectors of soft information about the individuals. That component about the synergies is defined from multiple vectors. The generation of these vectors is based on the use of  trapezoidal fuzzy sets \cite{fuzzyNumbers}.} 

\subsection{Experiment Design}\label{subsec:experiment}

 Following the idea in \cite{biological}, then, we explain how we generate the benchmark models. Each one will represent an MEFVFG with $n=256$ nodes. This process has two main parts: the definition of the adjacency matrix and the generation of the additional information.
 
 To approach the manipulation of multiple vectors from a benchmarking perspective, we propose the following: in every vector, the value of each component   depends on certain trapezoidal fuzzy sets, specifically saying \textit{low} and \textit{high} fuzzy sets. \textit{Low} fuzzy sets are related to the generation of the components of each vector which imply scarce connections among the nodes, whereas \textit{high} fuzzy sets refer to the generation of the components of the vectors which  imply many connections among the nodes. Therefore,  in each vector, the component related to nodes which are in the same community are generated as \textit{high} fuzzy sets, whereas the components related to nodes of different communities are generated as \textit{low} fuzzy sets. Let us emphasize that in the simulation process presented, what we randomly generate are the values $D(\widetilde{f}^{i})$ as \textit{high} or \textit{low} depending on the trapezoidal fuzzy sets $\widetilde{f}^{i}$, being $D$ a defuzzification operator.  
 
 We have $r$ vectors of fuzzy sets as the starting point in each benchmark, where $r$ is the amount of communities embedded in the synergies matrix, $F$.  Each vector is associated with a community $C_i$, so nodes belonging to $C_i$  have a \textit{high} value in the vector $\widetilde{f}^i$, whereas nodes which are not in $C_i$  have a \textit{low} value in $\widetilde{f}^i$: $D(\widetilde{f}^i_{j})=\hbar$, if $j \in C_i$; $D(\widetilde{f}^i_{j})=\ell$, if $j \notin C_i$. The process for defining and simulating these trapezoidal fuzzy sets is detailed below. To analyze different scattering  of the $\ell$ and $\hbar$ fuzzy sets, several combinations of the parameters $a$, $b$, $c$ and $d$ are considered (see Figures \ref{cap3:fig:fuzzyBajo} and  \ref{cap3:fig:fuzzyAlto}). For example, to define a benchmark graph with $four$ communities, we have to generate $four$ $n$-vectors $\left(D(\widetilde{f}^1),D(\widetilde{f}^2),D(\widetilde{f}^3),D(\widetilde{f}^4)\right)$ with $n=256$ nodes. 
 \color{black}

{ Each benchmark model represents an MEFVFG \color{black} summarized into two matrices: one of direct connections (adjacency $A$) and another of additional information (synergies matrix $F$, obtained from the soft information vectors).} Below, we explain the generation of them.

\begin{enumerate}
\item \textbf{Adjacency matrix.} The adjacency matrix $A$ is randomly generated according to Equation \eqref{eq:benchmarkingGeneracion} for a set $V$ with $256$ nodes. We consider different combinations of the values of  parameters $\alpha$ and $\beta$ regarding the input/output values ($z_{in}$ and $z_{out}$), as shown in Table~\ref{cap3:ta:benchmark} (similarly to the proposal in \cite{biological}). These parameters regulate the density of the connections matrix, $A$, whose generation process is shown in Algorithm \ref{alg:generaAdy}.

\begin{equation}\label{eq:benchmarkingGeneracion}
P(i,j)= \left\lbrace 
\begin{array}{lcl} 
\alpha  & \textup{if } & i,j \in C_k\\ 
\beta & \textup{if } &  \textup{otherwise}
\end{array} \right.
\end{equation}

\begin{table}[H]
\caption{\hl{Parameters used} to generate the adjacency matrix $A$ of each model.}\label{cap3:ta:benchmark}
 \centering
  \begin{adjustwidth}{-\extralength}{0cm}
 \scalebox{1}{\begin{tabularx}{\fulllength}{cccccccccc} \toprule
& \textbf{\begin{tabular}[c]{@{}c@{}} Network 1\end{tabular}} &
\textbf{\begin{tabular}[c]{@{}c@{}}Network 2\end{tabular}} &
\textbf{\begin{tabular}[c]{@{}c@{}}Network 3\end{tabular}} &
\textbf{\begin{tabular}[c]{@{}c@{}}Network 4\end{tabular}} &
\textbf{\begin{tabular}[c]{@{}c@{}}Network 5\end{tabular}} &
\textbf{\begin{tabular}[c]{@{}c@{}}Network 6\end{tabular}} &
\textbf{\begin{tabular}[c]{@{}c@{}}Network 7\end{tabular}} &
\textbf{\begin{tabular}[c]{@{}c@{}}Network 8\end{tabular}} &
\textbf{\begin{tabular}[c]{@{}c@{}}Network 9\end{tabular}} \\\midrule
${\alpha}$ & 0.45 & 0.4 & 0.35 & 0.325 & 0.3 & 0.275 & 0.25 &  0.225 & 0.2 \\ 
${\beta}$ & 0.016 &  0.033 & 0.05 &  0.058 & 0.066 & 0.075 & 0.083 & 0.091 & 0.1 \\ \bottomrule
\end{tabularx}}\end{adjustwidth}
\end{table}
\vspace{-9pt}
\begin{algorithm}[H]
 \begin{algorithmic}[1]
  \STATE{\textbf{Input}: $\left(|C_1|, \dots |C_r|\right), \alpha, \beta, n$;}
  \STATE{\textbf{Output}: $A$;}
  \STATE{$A\left(i,j\right)\leftarrow 0$, $\forall i,j=1,\dots , n$;}
  \FOR{$(i=1)$ \TO $(n)$}
  \FOR{$(i=1)$ \TO $(n)$}
  \FOR{$(\ell=1)$ \TO $(r)$}
  \STATE{$\epsilon \leftarrow rand(0,1)$;}
  \IF{$\left(|C_{\ell-1}|<i\leq|C_{\ell}|\right)$ $and$ $\left(|C_{\ell-1}|<j\leq|C_{\ell}|\right)$}
  \IF{$\epsilon< \alpha$}
  \STATE{$A(i,j) \leftarrow 1;$}
  \ENDIF
  \ELSE
  \IF{$\epsilon< \beta$}
  \STATE{$A(i,j) \leftarrow  1;$}
  \ENDIF
  \ENDIF
  \ENDFOR
  \ENDFOR
  \ENDFOR
  \STATE\textbf{return}$(A)$;
 \end{algorithmic}\caption{\textit{Generate Adjacency}}
 \label{alg:generaAdy}
\end{algorithm}
\item {\textbf{{Low trapezoidal fuzzy sets generation}.} This type of fuzzy sets $\widetilde{f}^{i}$, shown in Figure \ref{cap3:fig:fuzzyBajo}, are generated to represent, in each vector $D(\widetilde{f}^{i})$, the components related to the elements with a \textit{low} value in the characteristic of the corresponding vector.}
\vspace{-6pt}
\begin{figure}[H]
\begin{minipage}{0.33\textwidth}
    \includegraphics[scale=1.5]{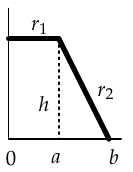}
    \end{minipage}
\caption{\textit{Low} trapezoidal fuzzy {set}.}\label{cap3:fig:fuzzyBajo}
\end{figure}


After fixing the values $a$ and $b$, we can calculate  the lines $r_1$ and $r_2$ to obtain a trapezoid below them with area $1$. 
  Particularly,  $r_1$ is defined as $y=h$, where $h$ is a value chosen so that the value of the corresponding integral is $1$:   $1=a\* h+\left(\frac{b-a}{2}\right)\*h \Longrightarrow h=\frac{2}{a+b}$. 
  
  On the other hand, $r_2$ is the line through the points $\left(\frac{2}{a+b},a\right)$ and $\left(0,b\right)$, so 
  
  \begin{equation*}
  r_2=\begin{cases}
  \frac{2}{a+b}=\alpha+\beta\*a \\
  0=\alpha+\beta\*b
  \end{cases}
  \end{equation*}
  
  By isolating $\alpha$ and $\beta$, $\alpha=\frac{-2b}{\left(a+b\right)\left(a-b\right)}$ and $\beta=\frac{2}{\left(a+b\right)\left(a-b\right)}$,
the distribution function of the \textit{low} trapezoidal fuzzy set is:  
  \begin{equation*}
  F\left(x\right)=\begin{cases}
  \frac{2x}{a+b}, & x\in[0,a] \\
  \frac{2a}{a+b}+\int_{a}^{x}\left(\frac{-2b}{\left(a+b\right)\left(a-b\right)}+\frac{2z}{\left(a+b\right)\left(a-b\right)}\right)dz, &  x\in(a,b] \end{cases}
  \end{equation*}
  
  where  
  $$\frac{2a}{a+b}+\int_{a}^{x}\left(\frac{-2b}{\left(a+b\right)\left(a-b\right)}+\frac{2z}{\left(a+b\right)\left(a-b\right)}\right)dz=\frac{2a}{a+b}+\frac{\left(x-b\right)^2-\left(a-b\right)^2}{\left(a+b\right)\left(a-b\right)}$$
 
 $ $
 
 Once the \textit{low} fuzzy {set} is characterized, in the following denoted by $\ell$, we apply the inverse method. First, we have to calculate the inverse function of $F$,  $F^{-1}\left(x\right)$; then, we simulate a value between $0$ and $1$ ($p$). Finally, $F^{-1}\left(p\right)$ is the value assigned to an edge which connect nodes which are not in the same community.
  
  \begin{itemize}
\item If $p\leq\frac{2a}{a+b}\Longrightarrow p=\frac{2x}{a+b} \Longrightarrow x=\frac{\left(a+b\right)p}{2}$.
\item If $p>\frac{2a}{a+b}\Longrightarrow p=\frac{2a}{a+b}+\frac{\left(x-b\right)^2-\left(a-b\right)^2}{\left(a+b\right)\left(a-b\right)} \\ \Longrightarrow x=b-\sqrt{\left(p-\frac{2a}{a+b}\right)\left(a+b\right)\left(a-b\right)+\left(a-b\right)^2}$. 

$ $

We take the sign `$-$' because $x-b<0$.
  \end{itemize}
 Then, if  $p\leadsto U(0,1)$, the \textit{low} values considered are obtained as:
 \begin{small}
  \begin{equation}\begin{cases}
  x=\frac{\left(a+b\right)p}{2} & \textup{ if } p\leq\frac{2a}{a+b} \\
  x=b-\sqrt{\left(p-\frac{2a}{a+b}\right)\left(a+b\right)\left(a-b\right)+\left(a-b\right)^2} & \textup{ otherwise}
  \end{cases}\end{equation}
  \end{small}

\noindent This process is summarized in Algorithm \ref{alg:fuzzyBajo}.
\begin{algorithm}[H]
 \begin{algorithmic}[1]
  \STATE\textbf{Input}: $a, b$;
  \STATE\textbf{Output}: $\ell$;
  \STATE{$p \leftarrow rand(0,1)$;}
  \IF{$p\leq \frac{2a}{a+ b}$}
  \STATE{$\ell \leftarrow \frac{\left(a+b\right)p}{2}$;}
  \ELSE
  \STATE{$\ell \leftarrow  b - \sqrt{\left(p-\frac{2a}{a+b}\right)\left(a+b\right)\left(a-b\right) + \left(a+b\right)^2}$;}
  \ENDIF
  \STATE\textbf{return}$\left(\ell\right)$;
 \end{algorithmic}\caption{\textit{{Low Fuzzy Set}}}
 \label{alg:fuzzyBajo}
\end{algorithm}

\item \textbf{High trapezoidal fuzzy sets generation.} {This type of fuzzy sets, $\widetilde{f}^{i}$, shown in Figure \ref{cap3:fig:fuzzyAlto}, are generated to represent, in each vector $D(\widetilde{f}^i)$, the components related to the elements with a \textit{high} value in the characteristic of the corresponding vector. }

\begin{figure}[H]
\begin{minipage}{0.33\textwidth}
    \includegraphics[scale=1.5]{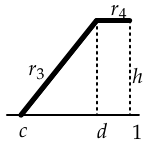}
    \end{minipage}
\caption{\textit{High} trapezoidal fuzzy {set}.}\label{cap3:fig:fuzzyAlto}
\end{figure}

  
After fixing the values $c$ and $d$, we can calculate  the lines $r_3$ and $r_4$ to obtain a trapezoid below them with area $1$. 
 Particularly, $r_3$ is defined as $y=h$, where $h$ is a value chosen so that the value of the corresponding integral is $1$:
  $1=\frac{\left(d-c\right)\*h}{2}\left(1-d\right)\*h \Longrightarrow h=\frac{2}{\left(1-d\right)+\left(1-c\right)}$
  
  On the other hand, $r_4$ is the line through points $\left(\frac{2}{\left(1-d\right)+\left(1-c\right)},d\right)$ and $\left(0,c\right)$, so:
  
  \begin{equation*}
  r_4=\begin{cases}
  \frac{2}{\left(1-d\right)+\left(1-c\right)}=\alpha+\beta\*d \\
  0=\alpha+\beta\*c
  \end{cases}
  \end{equation*}
  
  By isolating $\alpha$ and $\beta$, $\alpha=\frac{-2c}{\left[\left(1-d\right)+\left(1-c\right)\right]\left(d-c\right)}$ and $\beta=\frac{2}{\left[\left(1-d\right)+\left(1-c\right)\right]\left(d-c\right)}$, so the \linebreak distribution function of the \textit{high} trapezoidal fuzzy set is:
  \begin{small}
  \begin{equation*}
  F\left(x\right)=\begin{cases}
  \int_{c}^{x}\left(\frac{-2c+2z}{\left(\left(1-d\right)+\left(1-c\right)\right)\left(d-c\right)}\right)dz, & x\in[c,d] \\
  \int_{c}^{d}\left(\frac{-2c+2z}{\left(\left(1-d\right)+\left(1-c\right)\right)\left(d-c\right)}\right)dz+\int_{z=d}^{z=x}\frac{2}{\left(1-d\right)+\left(1-c\right)}, & x\in(d,1] \end{cases}
  \end{equation*}
  \end{small}
  where


    \begin{itemize}\item $\int_{c}^{x}\frac{-2c+2z}{\left(\left(1-d\right)+\left(1-c\right)\right)\left(d-c\right)}dz=\frac{\left(x-c\right)^2}{\left(\left(1-d\right)+\left(1-c\right)\right)\left(d-c\right)}$;
    \item $\int_{c}^{d}\left(\frac{-2c+2z}{\left(\left(1-d\right)+\left(1-c\right)\right)\left(d-c\right)}\right)dz+\int_{d}^{x}\frac{2}{\left(1-d\right)+\left(1-c\right)}=\frac{\left(x-d\right)+\left(x-c\right)}{\left(1-d\right)+\left(1-c\right)}$.\end{itemize}
    
  $ $
  
As with \textit{low} fuzzy {sets}, we apply the inverse method to simulate the values of the \textit{high}  fuzzy {sets} ($\hbar$ in the following). Then,  the value $F^{-1}\left(p\right)$ is:

$ $

\begin{itemize}
\item If $p\leq\frac{d-c}{\left(1-d\right)+\left(1-c\right)} \ \Longrightarrow \  p=\frac{\left(x-c\right)^2}{\left(\left(1-d\right)+\left(1-c\right)\right)\left(d-c\right)} \\ \Longrightarrow \ x=c+\sqrt{p\left(d-c\right)\left[\left(1-d\right)+\left(1-c\right)\right]}$;
\item If $p>\frac{d-c}{\left(1-d\right)+\left(1-c\right)}\ \Longrightarrow \  p=\frac{\left(x-d\right)+\left(x-c\right)}{\left(1-d\right)+\left(1-c\right)} \Longrightarrow \ x=\frac{p\left(\left(1-d\right)+\left(1-c\right)+d+c\right)}{2}$.

We take the sign `$+$' because $x-d>0$.
  \end{itemize} 
  
  $ $
  
Then, if  $p\leadsto  U(0,1)$, the \textit{high} values considered are obtained as:

\begin{small}
\begin{equation}\begin{cases}
x=c+\sqrt{p\left(d-c\right)\left(\left(1-c\right)+\left(1-d\right)\right)} & \textup{ if } p\leq\frac{d-c}{\left(1-c\right)+\left(1-d\right)} \\
x=\frac{p\left(\left(1-c\right)+\left(1-d\right)\right)+c+d}{2} & \textup{ otherwise}\end{cases}\end{equation}\end{small}
 
\noindent This process is summarized in Algorithm \ref{alg:fuzzyGrande}.

 \begin{algorithm}[H]
  \begin{algorithmic}[1]
\STATE\textbf{Input}: $c, d$;
\STATE\textbf{Output}: $\hbar$;
\STATE{$p \leftarrow rand(0,1)$;}
\IF{$\left(p\leq \frac{d-c}{(1-c)+(1-d)}\right)$}
\STATE{$\hbar \leftarrow c+\sqrt{p\*(d-c)\*((1-c)+(1-d))}$;}
\ELSE
\STATE{$\hbar \leftarrow \frac{p\left(\left(1-c\right)+\left(1-d\right)\right) +c+d}{2}$;}
\ENDIF
\STATE\textbf{return}$\left(\hbar\right)$;
  \end{algorithmic}\caption{\textit{High Fuzzy Set}}
  \label{alg:fuzzyGrande}
 \end{algorithm}

 \item \textbf{Generate multiple vectors.} 
 In each benchmark model, we have $r$ vectors as the starting point, where $r$ is the amount of communities embedded in the synergies matrix, $F$. Each vector is associated with a community $C_i$, so that nodes belonging to $C_i$ will have a \textit{high} value in $\widetilde{f}^i$, whereas the nodes which are not in $C_i$ will have a \textit{low} value in $\widetilde{f}^i$ {(then $D(\widetilde{f}^i_{j})=\hbar$, if $j \in C_i$; $D(\widetilde{f}^i_{j})=\ell$, if $j \notin C_i$).} Different combinations of the parameters $a$, $b$, $c$ and $d$ are considered to generate the \textit{low}/\textit{high} fuzzy sets (see Table \ref{cap3:ta:parametrosRel}). These combinations affect the scattering of the $\ell$ and $\hbar$ fuzzy sets.  The process is summarized in Algorithm \ref{alg:generaVec}.
\begin{algorithm}[H]
 \begin{algorithmic}[1]
  \STATE{\textbf{Input}: $\left(|C_1|, \dots |C_r|\right), a, b, c, d$;}
  \STATE{\textbf{Output}: $multipleVectors$;}
  \STATE{$|C_0|\leftarrow 0$;}
  \STATE {{$multipleVectors \leftarrow 0;$ (matrix $r\times n$, the line $\ell$ represents the vector $D(\widetilde{f}^\ell)$)}}
  \FOR{$(\ell=1)$ \TO $(r)$}  
  \FOR{$(i=1)$ \TO $(n)$}
  \IF{$|C_{\ell-1}|<i\leq|C_{\ell}|$}
  \STATE{$multipleVectors(\ell, i) \leftarrow High Fuzzy set(c,d)$;}
  \ELSE
  \STATE{$multipleVectors(\ell, i) \leftarrow Low Fuzzy set(a,b)$;}
  \ENDIF  
  \ENDFOR
  \ENDFOR
  \STATE\textbf{return}$(multipleVectors)$;
 \end{algorithmic}\caption{\textit{Generate Multiple Vectors}}
 \label{alg:generaVec}
\end{algorithm}


\vspace{-6pt}
\begin{table}[H]
\caption{Parameters to generate the matrix $F$ of the benchmark model.}
 \label{cap3:ta:parametrosRel} 
 \scalebox{1}{%
\begin{tabularx}{\textwidth}{cccccccccc}
\toprule
&
\textbf{\begin{tabular}[c]{@{}c@{}} Case 1\end{tabular}} &
\textbf{\begin{tabular}[c]{@{}c@{}}Case 2\end{tabular}} &
\textbf{\begin{tabular}[c]{@{}c@{}}Case 3\end{tabular}} &
\textbf{\begin{tabular}[c]{@{}c@{}}Case 4\end{tabular}} &
\textbf{\begin{tabular}[c]{@{}c@{}}Case 5\end{tabular}} &
\textbf{\begin{tabular}[c]{@{}c@{}}Case 6\end{tabular}} &
\textbf{\begin{tabular}[c]{@{}c@{}}Case 7\end{tabular}} &
\textbf{\begin{tabular}[c]{@{}c@{}}Case 8\end{tabular}} &
\textbf{\begin{tabular}[c]{@{}c@{}}Case 9\end{tabular}} \\ \midrule
\multicolumn{1}{c}{\textbf{$\mathbf{a}$}} & 0   & 0   & 0   & 0.1 & 0.1 & 0.1 & 0.2 & 0.2 & 0.2 \\ 
\multicolumn{1}{c}{\textbf{$\mathbf{b}$}} & 0.1 & 0.1 & 0.1 & 0.2 & 0.2 & 0.2 & 0.3 & 0.3 & 0.3 \\ 
\multicolumn{1}{c}{\textbf{$\mathbf{c}$}} & 0.9 & 0.8 & 0.7 & 0.9 & 0.8 & 0.7 & 0.9 & 0.8 & 0.7 \\ 
\multicolumn{1}{c}{\textbf{$\mathbf{d}$}} & 1   & 0.9 & 0.8 & 1   & 0.9 & 0.8 & 1   & 0.9 & 0.8 \\ \bottomrule 
  \end{tabularx}%
 }
\end{table}
\item \textbf{Synergies matrix.} {From  vectors generated with the Algorithm \textit{Generate Multiple Vectors}, we obtain  $\left(\mu_{f^1}^a, \dots , \mu_{f^r}^a\right)$.  
 We consider the matrices $\left({F}^{1}, \dots, {F}^{r}\right)$ and the adjacency of the corresponding MAWG.
The second component of each benchmark is an aggregation of these matrices, $F=\Phi\left({F}^{1}, \dots, {F}^{r}\right) = \max\left({F}^{1}, \dots, {F}^{r}\right)$.
We summarize this process in Algorithm \ref{alg:generaMatVector} for the particular case $p=1$}.

\begin{algorithm}[H]
 \begin{algorithmic}[1]
  \STATE{\textbf{Input}: $\left(|C_1|, \dots |C_r|\right), a, b, c, d$;}
  \STATE{\textbf{Output}: $F$;}
  \STATE{$multipleVectors \leftarrow Generate Multiple Vectors\left(\left(|C_1|, \dots |C_r|\right), a, b, c, d\right);$}
  \FOR{$(\ell=1)$ \TO $(r)$} 
  \FOR{$(i=1)$ \TO $(n)$}
  \STATE{$Sh(\ell,i) \leftarrow \frac{multipleVectors\left(\ell,i\right)}{\sum_{k=1}^n multipleVectors\left(\ell,k\right)};$ 
  }
  \FOR{$(j=1)$ \TO $(n)$}
  \STATE{$Sh_{j}(\ell,i) \leftarrow \frac{multipleVectors\left(\ell,i\right)}{\sum_{\substack{k \neq i\\k\in V}}^n multipleVectors\left(\ell,k\right)};$ 
  }
  \ENDFOR
  \ENDFOR
  \ENDFOR
  \FOR{$(\ell=1)$ \TO $(r)$} 
  \FOR{$(i=1)$ \TO $(n)$}
  \FOR{$(j=1)$ \TO $(n)$}
  \STATE{$F^{\ell}\left(i,j\right) \leftarrow \min\{|Sh(\ell,i)-Sh_{j}(\ell,i)|, |Sh(\ell,j)-Sh_{i}(\ell,j)|\}$};
  \ENDFOR
  \ENDFOR
  \ENDFOR
  \STATE $F \leftarrow \max\{F^1, \dots ,F^r\};$
  \STATE\textbf{return}$(F)$;
 \end{algorithmic}\caption{\textit{Matrix From Multiple Vectors}}
 \label{alg:generaMatVector}
\end{algorithm}
\end{enumerate}

\subsection{Results}\label{subsec:results}

Then, we show the evaluation of { the proposed methodology in the $1$-additive stage. We do this to avoid exponential complexity in computing fuzzy measures. Nevertheless, there is no reason to think that the goodness of the partitions obtained, and therefore the accuracy of the evaluated method, will worsen if non-additive measures are considered.} To do so, we consider several structures which vary in size and number of groups. Each of them represents an MEFVFG, $\widehat{G}=\left(V,E,\left(\mu_{f^1,p^1},\dots,\mu_{f^r,p^r}\right)\right)$ with two independent components. One of them, $G=\left(V,E\right)$, is related to the direct connections among the nodes represented by edges. The other,   $\left(\mu_{f^1,p^1},\dots,\mu_{f^r,p^r}\right)$, is used to define a relations matrix $F$.

For each combination of  $\alpha$ and $\beta$; $a$, $b$, $c$ and $d$,  we analyze the linear combination ${M}=\theta\left(A,F\right)=\gamma A + \left(1-\gamma\right)F$, by considering $\gamma=0$ (this is the only case in which, including the additional information, we can assert the partition which should be obtained).

In Tables  \ref{cap3:ta:modelo1}--\ref{cap3:ta:modelo4}, we show the average of the NMI obtained from $100$ iterations of each combination of $\alpha$ and $\beta$, concerning matrix $A$, and the parameters $a$, $b$, $c$ and $d$ for the definition of the vectors which give rise to the synergies matrix $F$. To simplify the interpretation of the results, these tables display the values in different colors: the closer the value is to $1$ (i.e., the better the result), the lighter the color. 

\begin{itemize}
\item \textbf{{Benchmark graph. Model  1.}}\label{cap3:com:bench1}
 \noindent \textls[-15]{It is the simpler benchmark model, showed in the Figure \ref{cap3:fig:caso1}. The adjacency matrix has two communities with an expected size of $128$ each, being $<k>=128\*\alpha+128\*\beta$  the expected degree of each node. $F_1$ is obtained from  vectors} {$\left(D(\widetilde{f}^1), D(\widetilde{f}^2), D(\widetilde{f}^3), D(\widetilde{f}^4)\right)$}, so the  $256$ \textls[-5]{nodes are organized into four groups $C_1^{{F}}, \dots , C_4^{{F}}$ with expected size} \mbox{$|C_i^{{F}}|=64$}.
 In Table \ref{cap3:ta:modelo1}, we show the results. Note that the tested algorithm always recovers the standard partition, even when the networks are sparse. 
 
 \begin{figure}[H]
\begin{minipage}{0.33\textwidth}
    \includegraphics[scale=1.0]{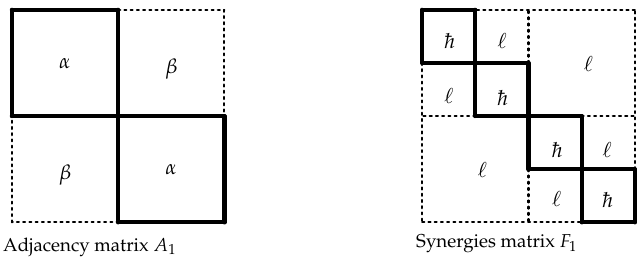}
    \end{minipage}
 \caption{Benchmark 
 graph. Model  1.}\label{cap3:fig:caso1}
\end{figure}

\vspace{-6pt}
 \begin{table}[H]
 \caption{\hl{NMI.} Model 1.} 
  \label{cap3:ta:modelo1}\centering
  \scalebox{1}{%
 \begin{tabularx}{\textwidth}{cccccccccc}
\toprule
\textbf{NMI} &
\textbf{\begin{tabular}[c]{@{}c@{}}$\mathbf{F_1}$\\ Case 1\end{tabular}} &
\textbf{\begin{tabular}[c]{@{}c@{}}$\mathbf{F_1}$\\ Case 2\end{tabular}} &
\textbf{\begin{tabular}[c]{@{}c@{}}$\mathbf{F_1}$\\ Case 3\end{tabular}} &
\textbf{\begin{tabular}[c]{@{}c@{}}$\mathbf{F_1}$\\ Case 4\end{tabular}} &
\textbf{\begin{tabular}[c]{@{}c@{}}$\mathbf{F_1}$\\ Case 5\end{tabular}} &
\textbf{\begin{tabular}[c]{@{}c@{}}$\mathbf{F_1}$\\ Case 6\end{tabular}} &
\textbf{\begin{tabular}[c]{@{}c@{}}$\mathbf{F_1}$\\ Case 7\end{tabular}} &
\textbf{\begin{tabular}[c]{@{}c@{}}$\mathbf{F_1}$\\ Case 8\end{tabular}} &
\textbf{\begin{tabular}[c]{@{}c@{}}$\mathbf{F_1}$\\ Case 9\end{tabular}} \\ \midrule
\textbf{\begin{tabular}[c]{@{}c@{}}$\mathbf{A_1}$\\ Network 1\end{tabular}} & \textcolor{grey5}{1 
} & \textcolor{grey5}{1} & \textcolor{grey5}{1} & \textcolor{grey5}{1} & \textcolor{grey5}{1} & \textcolor{grey5}{1} & \textcolor{grey5}{1} & \textcolor{grey5}{1} & \textcolor{grey5}{1} \\ 
\textbf{\begin{tabular}[c]{@{}c@{}}$\mathbf{A_1}$\\ Network 2\end{tabular}} & \textcolor{grey5}{1} & \textcolor{grey5}{1} & \textcolor{grey5}{1} & \textcolor{grey5}{1} & \textcolor{grey5}{1} & \textcolor{grey5}{1} & \textcolor{grey5}{1} & \textcolor{grey5}{1} & \textcolor{grey5}{1} \\ 
\textbf{\begin{tabular}[c]{@{}c@{}}$\mathbf{A_1}$\\ Network 3\end{tabular}} & \textcolor{grey5}{1} & \textcolor{grey5}{1} & \textcolor{grey5}{1} & \textcolor{grey5}{1} & \textcolor{grey5}{1} & \textcolor{grey5}{1} & \textcolor{grey5}{1} & \textcolor{grey5}{1} & \textcolor{grey5}{1} \\ 
\textbf{\begin{tabular}[c]{@{}c@{}}$\mathbf{A_1}$\\ Network 4\end{tabular}} & \textcolor{grey5}{1} & \textcolor{grey5}{1} & \textcolor{grey5}{1} & \textcolor{grey5}{1} & \textcolor{grey5}{1} & \textcolor{grey5}{1} & \textcolor{grey5}{1} & \textcolor{grey5}{1} & \textcolor{grey5}{1} \\ 
\textbf{\begin{tabular}[c]{@{}c@{}}$\mathbf{A_1}$\\ Network 5\end{tabular}} & \textcolor{grey5}{1} & \textcolor{grey5}{1} & \textcolor{grey5}{1} & \textcolor{grey5}{1} & \textcolor{grey5}{1} & \textcolor{grey5}{1} & \textcolor{grey5}{1} & \textcolor{grey5}{1} & \textcolor{grey5}{1} \\ 
\textbf{\begin{tabular}[c]{@{}c@{}}$\mathbf{A_1}$\\ Network 6\end{tabular}} & \textcolor{grey5}{1} & \textcolor{grey5}{1} & \textcolor{grey5}{1} & \textcolor{grey5}{1} & \textcolor{grey5}{1} & \textcolor{grey5}{1} & \textcolor{grey5}{1} & \textcolor{grey5}{1} & \textcolor{grey5}{1} \\ 
\textbf{\begin{tabular}[c]{@{}c@{}}$\mathbf{A_1}$\\ Network 7\end{tabular}} & \textcolor{grey5}{1} & \textcolor{grey5}{1} & \textcolor{grey5}{1} & \textcolor{grey5}{1} & \textcolor{grey5}{1} & \textcolor{grey5}{1} & \textcolor{grey5}{1} & \textcolor{grey5}{1} & \textcolor{grey5}{1} \\ 
\textbf{\begin{tabular}[c]{@{}c@{}}$\mathbf{A_1}$\\ Network 8\end{tabular}} & \textcolor{grey5}{1} & \textcolor{grey5}{1} & \textcolor{grey5}{1} & \textcolor{grey5}{1} & \textcolor{grey5}{1} & \textcolor{grey5}{1} & \textcolor{grey5}{1} & \textcolor{grey5}{1} & \textcolor{grey5}{1} \\ 
\textbf{\begin{tabular}[c]{@{}c@{}}$\mathbf{A_1}$\\ Network 9\end{tabular}} & \textcolor{grey5}{1} & \textcolor{grey5}{1} & \textcolor{grey5}{1} & \textcolor{grey5}{1} & \textcolor{grey5}{1} & \textcolor{grey5}{1} & \textcolor{grey5}{1} & \textcolor{grey5}{1} & \textcolor{grey5}{1} \\ \bottomrule
 \end{tabularx}%
  }
 \end{table}

\item \textbf{{Benchmark graph. Model  2.}}\label{cap3:com:bench2}
 \noindent Due to modularity resolution limit \cite{fortunatoResolutionLimit}, the Louvain algorithm is very sensitive to changes in the groups' size, particularly with small communities \cite{evolutionary}. Then, we test the algorithm in this context, where the groups of the standard structure are smaller than in Model 1, as it can be seen in the Figure \ref{cap3:fig:caso2}.   $F_2$ has eight communities  with $32$~nodes each.
 We also reduce  the communities in  $A_2$: it has four communities with $64$~nodes each, being $<k>=64\*\alpha+192\*\beta$ the expected degree.
 The obtained results are shown in Table \ref{cap3:ta:modelo2}. Despite the size reduction, our method provides good results.
 
  \begin{figure}[H]
\begin{minipage}{0.33\textwidth}
    \includegraphics[scale=1.0]{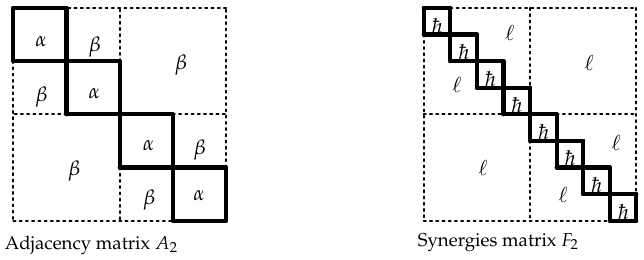}
    \end{minipage}
 \caption{Benchmark 
 graph. Model  2.}\label{cap3:fig:caso2}
\end{figure}

 \begin{table}[H]
 \caption{\hl{NMI.} Model 2.} 
  \label{cap3:ta:modelo2}
  \centering
  \scalebox{1}{%
 \begin{tabularx}{\textwidth}{c  c c c c c c c c c}
  \toprule
  \textbf{NMI} &
  \textbf{\begin{tabular}[c]{@{}c@{}}$\mathbf{F_2}$\\ Case 1\end{tabular}} &
  \textbf{\begin{tabular}[c]{@{}c@{}}$\mathbf{F_2}$\\ Case 2\end{tabular}} &
  \textbf{\begin{tabular}[c]{@{}c@{}}$\mathbf{F_2}$\\ Case 3\end{tabular}} &
  \textbf{\begin{tabular}[c]{@{}c@{}}$\mathbf{F_2}$\\ Case 4\end{tabular}} &
  \textbf{\begin{tabular}[c]{@{}c@{}}$\mathbf{F_2}$\\ Case 5\end{tabular}} &
  \textbf{\begin{tabular}[c]{@{}c@{}}$\mathbf{F_2}$\\ Case 6\end{tabular}} &
  \textbf{\begin{tabular}[c]{@{}c@{}}$\mathbf{F_2}$\\ Case 7\end{tabular}} &
  \textbf{\begin{tabular}[c]{@{}c@{}}$\mathbf{F_2}$\\ Case 8\end{tabular}} &
  \textbf{\begin{tabular}[c]{@{}c@{}}$\mathbf{F_2}$\\ Case 9\end{tabular}} \\ \midrule
  \textbf{\begin{tabular}[c]{@{}c@{}}$\mathbf{A_2}$\\ Network 1\end{tabular}} &  \textcolor{grey3}{1}  &  \textcolor{grey3}{1}  &  \textcolor{grey3}{1}  &  \textcolor{grey3}{1}  &  \textcolor{grey3}{1}  &  \textcolor{grey3}{1}  & \textcolor{grey8}{0.9987} & 0.9981 & \textcolor{grey5}{0.9994} \\ 
  \textbf{\begin{tabular}[c]{@{}c@{}}$\mathbf{A_2}$\\ Network 2\end{tabular}} &  \textcolor{grey3}{1}  &  \textcolor{grey3}{1}  &  \textcolor{grey3}{1}  &  \textcolor{grey3}{1}  &  \textcolor{grey3}{1}  &  \textcolor{grey3}{1}  & \textcolor{grey8}{0.9986} & 0.9980  & \textcolor{grey5}{0.9994} \\ 
  \textbf{\begin{tabular}[c]{@{}c@{}}$\mathbf{A_2}$\\ Network 3\end{tabular}} &  \textcolor{grey3}{1}  &  \textcolor{grey3}{1}  &  \textcolor{grey3}{1}  &  \textcolor{grey3}{1}  &  \textcolor{grey3}{1}  &  \textcolor{grey3}{1}  & \textcolor{grey5}{0.9993} & \textcolor{grey5}{0.9992} & \textcolor{grey5}{0.9994} \\ 
  \textbf{\begin{tabular}[c]{@{}c@{}}$\mathbf{A_2}$\\ Network 4\end{tabular}} &  \textcolor{grey3}{1}  &  \textcolor{grey3}{1}  &  \textcolor{grey3}{1}  &  \textcolor{grey3}{1}  &  \textcolor{grey3}{1}  &  \textcolor{grey3}{1}  & \textcolor{grey8}{0.9986} & 0.9980  & \textcolor{grey5}{0.9996} \\ 
  \textbf{\begin{tabular}[c]{@{}c@{}}$\mathbf{A_2}$\\ Network 5\end{tabular}} &  \textcolor{grey3}{1}  &  \textcolor{grey3}{1}  &  \textcolor{grey3}{1}  &  \textcolor{grey3}{1}  &  \textcolor{grey3}{1}  &  \textcolor{grey3}{1}  & 0.9984 & \textcolor{grey5}{0.9990}  & \textcolor{grey5}{0.9991} \\ 
  \textbf{\begin{tabular}[c]{@{}c@{}}$\mathbf{A_2}$\\ Network 6\end{tabular}} &  \textcolor{grey3}{1}  &  \textcolor{grey3}{1}  &  \textcolor{grey3}{1}  &  \textcolor{grey3}{1}  &  \textcolor{grey3}{1}  &  \textcolor{grey3}{1}  & \textcolor{grey8}{0.9986} & 0.9984 & \textcolor{grey5}{0.9990}  \\ 
  \textbf{\begin{tabular}[c]{@{}c@{}}$\mathbf{A_2}$\\ Network 7\end{tabular}} &  \textcolor{grey3}{1}  &  \textcolor{grey3}{1}  &  \textcolor{grey3}{1}  &  \textcolor{grey3}{1}  &  \textcolor{grey3}{1}  &  \textcolor{grey3}{1}  & \textcolor{grey5}{0.9990}  & \textcolor{grey5}{0.9992} & \textcolor{grey5}{0.9992} \\ 
  \textbf{\begin{tabular}[c]{@{}c@{}}$\mathbf{A_2}$\\ Network 8\end{tabular}} &  \textcolor{grey3}{1}  &  \textcolor{grey3}{1}  &  \textcolor{grey3}{1}  &  \textcolor{grey3}{1}  &  \textcolor{grey3}{1}  &  \textcolor{grey3}{1}  & \textbf{0.9968} & \textcolor{grey5}{0.9992} & \textcolor{grey8}{0.9989} \\ 
  \textbf{\begin{tabular}[c]{@{}c@{}}$\mathbf{A_2}$\\ Network 9\end{tabular}} &  \textcolor{grey3}{1}  &  \textcolor{grey3}{1}  &  \textcolor{grey3}{1}  &  \textcolor{grey3}{1}  &  \textcolor{grey3}{1}  &  \textcolor{grey3}{1}  & \textcolor{grey5}{0.9993} & \textcolor{grey5}{0.9992} & \textcolor{grey5}{0.9996} \\ \bottomrule
 \end{tabularx}%
  }
 \end{table} 
 
\item \textbf{{Benchmark graph. Model 3.}}\label{cap3:com:bench3}
 \noindent Previous results set light on the high quality of the tested method in symmetric structures. However, the interest of every method goes further than synthetic structures; the main objective is to reach proper results in real cases. Then, we work with asymmetric structures to simulate more realistic networks, as it can be seen in the Figure \ref{cap3:fig:caso3}.  $F_3$ has four communities whose sizes are $|C_1^{{F}}|=43$, $|C_2^{{F}}|=42$, $|C_3^{{F}}|=43$, $|C_4^{{F}}|=96$, $|C_5^{{F}}|=32$. On the other hand, $A_3=A_1$. 
We show the results in Table \ref{cap3:ta:modelo3}.

\vspace{-6pt}

\begin{figure}[H]
\begin{minipage}{0.33\textwidth}
    \includegraphics[scale=1.0]{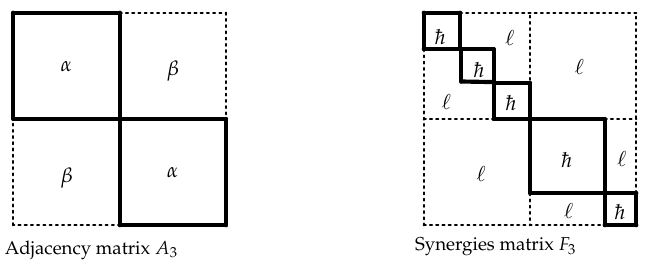}
    \end{minipage} 
     \caption{Benchmark 
 graph. Model  3.}\label{cap3:fig:caso3}
\end{figure}

\vspace{-6pt}
 \begin{table}[H]
 \caption{\hl{NMI.} Model 3.} 
  \label{cap3:ta:modelo3} 
  \scalebox{1}{%
 \begin{tabularx}{\textwidth}{c c c c c c c c c c}
 \toprule
 \textbf{NMI} &
 \textbf{\begin{tabular}[c]{@{}c@{}}$\mathbf{F_3}$\\ Case 1\end{tabular}} &
 \textbf{\begin{tabular}[c]{@{}c@{}}$\mathbf{F_3}$\\ Case 2\end{tabular}} &
 \textbf{\begin{tabular}[c]{@{}c@{}}$\mathbf{F_3}$\\ Case 3\end{tabular}} &
 \textbf{\begin{tabular}[c]{@{}c@{}}$\mathbf{F_3}$\\ Case 4\end{tabular}} &
 \textbf{\begin{tabular}[c]{@{}c@{}}$\mathbf{F_3}$\\ Case 5\end{tabular}} &
 \textbf{\begin{tabular}[c]{@{}c@{}}$\mathbf{F_3}$\\ Case 6\end{tabular}} &
 \textbf{\begin{tabular}[c]{@{}c@{}}$\mathbf{F_3}$\\ Case 7\end{tabular}} &
 \textbf{\begin{tabular}[c]{@{}c@{}}$\mathbf{F_3}$\\ Case 8\end{tabular}} &
 \textbf{\begin{tabular}[c]{@{}c@{}}$\mathbf{F_3}$\\ Case 9\end{tabular}} \\ \midrule
 \textbf{\begin{tabular}[c]{@{}c@{}}$\mathbf{A_3}$\\ Network 1\end{tabular}} &  \textcolor{grey5}{1}  &  \textcolor{grey5}{1}  &  \textcolor{grey5}{1}  &  \textcolor{grey5}{1}  &  \textcolor{grey5}{1}  &  \textcolor{grey5}{1}  & \textcolor{grey8}{0.9996} & 0.9994 & \textcolor{grey8}{0.9997} \\ 
 \textbf{\begin{tabular}[c]{@{}c@{}}$\mathbf{A_3}$\\ Network 2\end{tabular}} &  \textcolor{grey5}{1}  &  \textcolor{grey5}{1}  &  \textcolor{grey5}{1}  &  \textcolor{grey5}{1}  &  \textcolor{grey5}{1}  &  \textcolor{grey5}{1}  & \textcolor{grey8}{0.9997} &  \textcolor{grey5}{1}       &  \textcolor{grey5}{1}  \\ 
 \textbf{\begin{tabular}[c]{@{}c@{}}$\mathbf{A_3}$\\ Network 3\end{tabular}} &  \textcolor{grey5}{1}  &  \textcolor{grey5}{1}  &  \textcolor{grey5}{1}  &  \textcolor{grey5}{1}  &  \textcolor{grey5}{1}  &  \textcolor{grey5}{1}  & \textcolor{grey8}{0.9997} &  \textcolor{grey5}{1}       &  \textcolor{grey5}{1}       \\ 
 \textbf{\begin{tabular}[c]{@{}c@{}}$\mathbf{A_3}$\\ Network 4\end{tabular}} &  \textcolor{grey5}{1}  &  \textcolor{grey5}{1}  &  \textcolor{grey5}{1}  &  \textcolor{grey5}{1}  &  \textcolor{grey5}{1}  &  \textcolor{grey5}{1}  & \textbf{0.9990 
}  &  \textcolor{grey5}{1}       &  \textcolor{grey5}{1}       \\ 

\bottomrule
\end{tabularx}}
\end{table}

\begin{table}[H]\ContinuedFloat
\caption{{\em Cont.}}
 \label{cap3:ta:modelo3} 
  \scalebox{1}{%
 \begin{tabularx}{\textwidth}{c c c c c c c c c c}
 \toprule
 \textbf{NMI} &
 \textbf{\begin{tabular}[c]{@{}c@{}}$\mathbf{F_3}$\\ Case 1\end{tabular}} &
 \textbf{\begin{tabular}[c]{@{}c@{}}$\mathbf{F_3}$\\ Case 2\end{tabular}} &
 \textbf{\begin{tabular}[c]{@{}c@{}}$\mathbf{F_3}$\\ Case 3\end{tabular}} &
 \textbf{\begin{tabular}[c]{@{}c@{}}$\mathbf{F_3}$\\ Case 4\end{tabular}} &
 \textbf{\begin{tabular}[c]{@{}c@{}}$\mathbf{F_3}$\\ Case 5\end{tabular}} &
 \textbf{\begin{tabular}[c]{@{}c@{}}$\mathbf{F_3}$\\ Case 6\end{tabular}} &
 \textbf{\begin{tabular}[c]{@{}c@{}}$\mathbf{F_3}$\\ Case 7\end{tabular}} &
 \textbf{\begin{tabular}[c]{@{}c@{}}$\mathbf{F_3}$\\ Case 8\end{tabular}} &
 \textbf{\begin{tabular}[c]{@{}c@{}}$\mathbf{F_3}$\\ Case 9\end{tabular}} \\ \midrule
 \textbf{\begin{tabular}[c]{@{}c@{}}$\mathbf{A_3}$\\ Network 5\end{tabular}} &  \textcolor{grey5}{1}  &  \textcolor{grey5}{1}  &  \textcolor{grey5}{1}  &  \textcolor{grey5}{1}  &  \textcolor{grey5}{1}  &  \textcolor{grey5}{1}  &  \textcolor{grey5}{1}       & \textcolor{grey8}{0.9996} & \textbf{0.9992} \\ 
 \textbf{\begin{tabular}[c]{@{}c@{}}$\mathbf{A_3}$\\ Network 6\end{tabular}} &  \textcolor{grey5}{1}  &  \textcolor{grey5}{1}  &  \textcolor{grey5}{1}  &  \textcolor{grey5}{1}  &  \textcolor{grey5}{1}  &  \textcolor{grey5}{1}  &  \textcolor{grey5}{1}       & 0.9994 & 0.9994 \\ 
 \textbf{\begin{tabular}[c]{@{}c@{}}$\mathbf{A_3}$\\ Network 7\end{tabular}} &  \textcolor{grey5}{1}  &  \textcolor{grey5}{1}  &  \textcolor{grey5}{1}  &  \textcolor{grey5}{1}  &  \textcolor{grey5}{1}  &  \textcolor{grey5}{1}  & \textcolor{grey8}{0.9997} & \textcolor{grey8}{0.9996} & \textcolor{grey8}{0.9997} \\ 
 \textbf{\begin{tabular}[c]{@{}c@{}}$\mathbf{A_3}$\\ Network 8\end{tabular}} &  \textcolor{grey5}{1}  &  \textcolor{grey5}{1}  &  \textcolor{grey5}{1}  &  \textcolor{grey5}{1}  &  \textcolor{grey5}{1}  &  \textcolor{grey5}{1}  &  \textcolor{grey5}{1}       &  \textcolor{grey5}{1}       & 0.9994 \\ 
 \textbf{\begin{tabular}[c]{@{}c@{}}$\mathbf{A_3}$\\ Network 9\end{tabular}} &  \textcolor{grey5}{1}  &  \textcolor{grey5}{1}  &  \textcolor{grey5}{1}  &  \textcolor{grey5}{1}  &  \textcolor{grey5}{1}  &  \textcolor{grey5}{1}  &  \textcolor{grey5}{1}       &  \textcolor{grey5}{1}       & 0.9995 \\ \bottomrule
 \end{tabularx}%
  }
 \end{table} 
\item \textbf{{Benchmark graph. Model  4.}}\label{cap3:com:bench4}
 \noindent \textls[-20]{This model combines the reduction of the size} communities with partition asymmetry, as it can be seen in the Figure \ref{cap3:fig:caso4}. In this case,  $A_4=A_2$, and  $F_4$ has eight communities whose expected sizes are $|C_1^{{F}}|=24$, $|C_2^{{F}}|=40$, $|C_3^{{F}}|=64$, $|C_4^{{F}}|=21$, $|C_5^{{F}}|=22$, $|C_6^{{F}}|=21$, $|C_7^{{F}}|=32$ y $|C_8^{{F}}|=32$. Despite the obvious complexity of this structure, the results presented in Table \ref{cap3:ta:modelo4} show the good performance of the tested algorithm.
 
 \begin{figure}[H]
\begin{minipage}{0.33\textwidth}
    \includegraphics[scale=1.0]{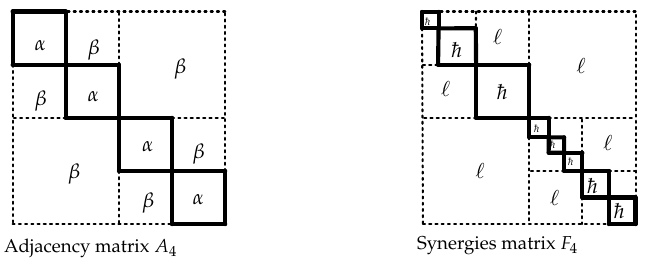}
    \end{minipage} 
     \caption{Benchmark 
 graph. Model  4.}\label{cap3:fig:caso4}
\end{figure}

 \vspace{-6pt}
 \begin{table}[H]
 \caption{\hl{NMI.} Model 4.}
  \label{cap3:ta:modelo4}
  \centering
  \scalebox{1}{%
 \begin{tabularx}{\textwidth}{c  c c c c c c c c c}
 \toprule
 \textbf{NMI} &
 \textbf{\begin{tabular}[c]{@{}c@{}}$\mathbf{F_4}$\\ Case 1\end{tabular}} &
 \textbf{\begin{tabular}[c]{@{}c@{}}$\mathbf{F_4}$\\ Case 2\end{tabular}} &
 \textbf{\begin{tabular}[c]{@{}c@{}}$\mathbf{F_4}$\\ Case 3\end{tabular}} &
 \textbf{\begin{tabular}[c]{@{}c@{}}$\mathbf{F_4}$\\ Case 4\end{tabular}} &
 \textbf{\begin{tabular}[c]{@{}c@{}}$\mathbf{F_4}$\\ Case 5\end{tabular}} &
 \textbf{\begin{tabular}[c]{@{}c@{}}$\mathbf{F_4}$\\ Case 6\end{tabular}} &
 \textbf{\begin{tabular}[c]{@{}c@{}}$\mathbf{F_4}$\\ Case 7\end{tabular}} &
 \textbf{\begin{tabular}[c]{@{}c@{}}$\mathbf{F_4}$\\ Case 8\end{tabular}} &
 \textbf{\begin{tabular}[c]{@{}c@{}}$\mathbf{F_4}$\\ Case 9\end{tabular}} \\ \midrule
 \textbf{\begin{tabular}[c]{@{}c@{}}$\mathbf{A_4}$\\ Network 1\end{tabular}} & \textcolor{grey8}{0.9975} & 0.9969 & \textcolor{grey8}{0.9974} & \textcolor{grey8}{0.9976} & \textcolor{grey8}{0.9970}  & 0.9968 & \textbf{0.9937} &  \textcolor{grey3}{1}       & \textcolor{grey5}{0.9994} \\ 
 \textbf{\begin{tabular}[c]{@{}c@{}}$\mathbf{A_4}$\\ Network 2\end{tabular}} & 0.9956 & 0.9960  & 0.9962 & \textcolor{grey8}{0.9977} & \textcolor{grey8}{0.9982} & \textcolor{grey8}{0.9981} & 0.9960  & \textcolor{grey5}{0.9991} &  \textcolor{grey3}{1}       \\ 
 \textbf{\begin{tabular}[c]{@{}c@{}}$\mathbf{A_4}$\\ Network 3\end{tabular}} & 0.9951 & 0.9964 & 0.9956 & \textcolor{grey8}{0.9980}  & \textcolor{grey8}{0.9972} & \textcolor{grey8}{0.9970}  & \textcolor{grey8}{0.9974} & \textcolor{grey8}{0.9975} & \textcolor{grey5}{0.9997} \\ 
 \textbf{\begin{tabular}[c]{@{}c@{}}$\mathbf{A_4}$\\ Network 4\end{tabular}} & \textbf{0.9949} & 0.9964 & 0.9963 & \textcolor{grey8}{0.9978} & 0.9963 & \textcolor{grey8}{0.9980}  & \textbf{0.9943} & \textcolor{grey5}{0.9990}  & \textcolor{grey8}{0.9973} \\ 
 \textbf{\begin{tabular}[c]{@{}c@{}}$\mathbf{A_4}$\\ Network 5\end{tabular}} & 0.9952 & 0.9956 & \textcolor{grey8}{0.9975} & 0.9960  & \textcolor{grey8}{0.9978} & \textcolor{grey8}{0.9976} & 0.9957 & \textcolor{grey8}{0.9988} & \textcolor{grey8}{0.9979} \\ 
 \textbf{\begin{tabular}[c]{@{}c@{}}$\mathbf{A_4}$\\ Network 6\end{tabular}} & 0.9954 & 0.9957 & \textcolor{grey8}{0.9972} & 0.9967 & 0.9954 & 0.9989 & \textcolor{grey8}{0.9973} & 0.9961 & \textcolor{grey5}{0.9990}  \\ 
 \textbf{\begin{tabular}[c]{@{}c@{}}$\mathbf{A_4}$\\ Network 7\end{tabular}} & \textbf{0.9944} & 0.9955 & \textcolor{grey8}{0.9978} & \textcolor{grey8}{0.9971} & 0.9969 & \textcolor{grey8}{0.9971} & \textcolor{grey8}{0.9980}  & \textcolor{grey8}{0.9989} &  \textcolor{grey3}{1}       \\ 
 \textbf{\begin{tabular}[c]{@{}c@{}}$\mathbf{A_4}$\\ Network 8\end{tabular}} & 0.9952 & 0.9958 & \textcolor{grey8}{0.9973} & \textcolor{grey8}{0.9974} & 0.9969 & \textcolor{grey8}{0.9980}  & \textcolor{grey5}{0.9996} & \textcolor{grey8}{0.9984} &  \textcolor{grey3}{1}       \\ 
 \textbf{\begin{tabular}[c]{@{}c@{}}$\mathbf{A_4}$\\ Network 9\end{tabular}} & 0.9952 & 0.9967 & \textcolor{grey8}{0.9974} & 0.9965 & 0.9969 & \textcolor{grey8}{0.9972} & \textcolor{grey8}{0.9972} & 0.9968 &  \textcolor{grey3}{1}       \\ \bottomrule
 \end{tabularx}%
  }  
 \end{table} 
 
 \end{itemize}
 
\section{Discussion and Conclusions}\label{sec:conclu} 

Within the framework of the analysis of networks and social networks, this work is based on the analysis of fuzzy information. We can find different main contributions in this article. Our first objective is the definition of a new representation model, whose base is made up of two types of information sources. The first is a set of individuals whose direct connections are known and represented by a crisp graph or network. The other part deals with some additional knowledge about these individuals in terms of soft information represented by family of vectors of fuzzy sets. On this basis, we define the multi-dimensional extended fuzzy graph based on fuzzy vectors (MEFVFG). This new model combines the crisp information of a graph with the soft information provided by the vectors about the individuals. We define it from the simplest case where there is only one vector of fuzzy sets to the multidimensional scenario involving multiple vectors.

Another goal of this work is the proposal of a specific application of this new model related to community detection. This problem has been previously addressed in terms of fuzzy measures that provide additional information about the synergies among the individuals in some works \cite{ampliinfus,multipleBip,polarizacion:covid}. In \cite{infus}, we proposed a methodology, named \textit{Duo-Louvain}, to find a ``good'' partition of the individuals of an extended fuzzy graph considering both the connections defined by the edges and also the additional information provided by the fuzzy measures. It was based on the well-known Louvain method \cite{blondel}, a greedy multi-phase algorithm based on local moving \cite{smart} and modularity optimization~\cite{girvanNewman}. That proposal \cite{infus} is the inspiration of this paper. Now, we work on the community detection in networks by considering additional soft information about the individuals of the network. Specifically, we face the existence of several fuzzy sets related to the nodes of a network, so our proposed application of the MEFVFG in current work is to obtain realistic communities in it. That idea is quite useful and goes beyond any other previous proposal, as it can be applied in a wide range of scenarios, for example, when any linguistic variable(s) appear. As far as we know, this situation has never been faced before, so this work intrinsically leads to the definition of a new type of problem.

Another important objective is the evaluation of the developed methodology. As mentioned above, the problem presented in this article has not been addressed before in the literature, so there are no other methods with which we can compare our proposal. Then, to evaluate the performance of the new algorithm, we present some experimental results developed on the basis of benchmarking \cite{benchmarking2} and NMI calculation \cite{nmi}. We develop some methods based on trapezoidal fuzzy sets with which we generate the elements of the gold models considered, each with a standard partition summarizing an MEFVFG, which should be detected by the evaluated algorithm. The high level of the results shown in  Section \ref{sec:computational}  allows us to assert the good performance of the proposed method:  the NMI calculated in almost all the scenarios is $1$, which means that our algorithm perfectly detects the standard partition despite the complexity of the considered model.

As further research, we stress the importance of an in-depth analysis of the distance between fuzzy sets. Specifically, we are interested in analyzing how far two fuzzy sets are in order to compute new measures of additional information to be later considered in the presented methodology. Our idea is to work with the Hausdorff distance between two fuzzy sets \cite{fuzzyDistance}, based on the classic metric with the same name \cite{variational}, which is used in mathematics to quantify how far two subsets of a metric space are from each other. To approach this theoretical approach, it is essential to be familiar with the properties of fuzzy sets and also with various topological concepts related to the measurement of distances in different~spaces.

Another important line of future work is not so theoretical but applied. Let us emphasize the importance of applying this methodology in real-life cases in order to obtain realistic groups of individuals which not only consider the direct connections between them but also some additional soft information. Real problems  are too complex to be represented by a crisp graph alone.  The need to include as many sources of information as possible is clear. For example,  in the behavior of people, things are not ``black or white''. To the question ``do you agree with this law?'', the answer may be something like ``well, more or less but not quite''. 
The more capable we are of representing these situations in a model, the more realistic the results obtained will be.  We have to be prepared to understand, model, and analyze the fuzzy knowledge of real life, for example, by the consideration of linguistic terms. Undeniably, it is worth an in-depth analysis of the linguistic terms that accompany any real problem, for whose study the tools and methodology proposed here can be crucial. When facing this type of problem, it is vitally important to take into account its difficulty, both computationally and in terms of understanding. When fuzzy elements appear, it is essential to be well prepared to consider tools that mitigate the intrinsic difficulty.

\vspace{12pt}
\authorcontributions{Conceptualization, I.G., D.G. and J.C.; methodology, I.G., D.G. and J.C.; software, I.G. and J.C.; validation, D.G., J.C. and R.E.; formal analysis, D.G.; investigation,  I.G., D.G., J.C. and R.E.; resources,  D.G., J.C. and R.E.; data curation, I.G.; writing---original draft preparation, I.G., D.G. and J.C.; writing---review and editing, I.G., D.G., J.C. and R.E.; visualization, I.G.; supervision, D.G., J.C. and R.E.; project administration, D.G.; funding acquisition, D.G., J.C. and R.E. All authors have read and agreed to the published version of the manuscript.}

\funding {\hl{This research has been partially supported by the Government of Spain, Grant Plan Nacional de I+D+i, PID2020-116884GB-I00, PGC2018096509-B-I00.}}

\institutionalreview {\hl{Not applicable}}

\informedconsent {\hl{Not applicable}}

\dataavailability{\hl{Not applicable}}

\acknowledgments{\hl{This research has been partially supported by the Government of Spain, Grant Plan Nacional de I+D+i, PID2020-116884GB-I00.}}

\conflictsofinterest{The authors declare no conflict of interest. The funders had no role in the design of the study; in the collection, analyses, or interpretation of data; in the writing of the manuscript; or in the decision to publish the~results.} 
    
\begin{adjustwidth}{-\extralength}{0cm}
\reftitle{References}

\end{adjustwidth}
\end{document}